\documentclass{scrartcl}

\usepackage{amsmath,amssymb,amsthm}
\usepackage[english]{babel}
\usepackage[utf8x]{inputenc}

\usepackage{xcolor}
\usepackage{pgf}
\usepackage{tikz}
\usepackage{tabulary}
\usepackage[shortlabels]{enumitem}
\usepackage[framemethod=tikz]{mdframed}
\usepackage{caption}
\usepackage[section]{placeins}
\usepackage{pdfpages}
\usepackage{mathtools}

\usepackage{pgfplots}
\usepackage{pgfmath}
\pgfplotsset{compat=1.13}

\definecolor{GarnetRed}{RGB}{115,0,10}
\definecolor{BoxBg}{gray}{0.95}

\usepackage[colorlinks=false,linkcolor=blue,citecolor=blue]{hyperref}
\usepackage[capitalise,noabbrev]{cleveref}
\usepackage{bbm}

\usetikzlibrary{shapes.misc}

\tikzset{cross/.style={cross out, draw=black, minimum size=2*(#1-\pgflinewidth), inner sep=0pt, outer sep=0pt}, cross/.default={2.6pt}}

\newcommand{\N}{\mathbb{N}}

\newcommand{\Z}{\mathbb{Z}}

\newcommand{\R}{\mathbb{R}}
\newcommand{\C}{\mathbb{C}}

\newcommand{\abs}[1]{\left| #1 \right|}
\newcommand{\norm}[1]{\left|\hspace*{-0.1em}\left| #1 \right|\hspace*{-0.1em}\right|}

\newcommand{\bref}[1]{(\ref{#1})}

\newcommand{\dz}{\,\mathrm{d}z}

\newcommand{\dk}{\,\mathrm{d}k}

\newcommand{\ii}{\mathrm{i}}

\DeclarePairedDelimiter\floor{\lfloor}{\rfloor}

\newtheorem{theorem}{Theorem}[section]

\newtheorem{task}[theorem]{Task}

\theoremstyle{definition}

\graphicspath{{./Images/}}

\title{Lattice Multislice Algorithm for Fast Simulation of Scanning Transmission Electron Microscopy Images}
\author{Christian Doberstein and Peter Binev \thanks{%
   This research was supported in part by the ARO grant W911NF-20-1-0318 and the NSF Grants DMS 2038080, DMS 1925919, and DMS 2245097.
   }}
\begin{document}

\maketitle

\begin{abstract}
\textbf{Abstract:}
We introduce a new approach to the numerical simulation of the Scanning Transmission Electron
Microscopy images. The Lattice Multislice Algorithm (LMA) takes advantage of the fact that
the electron waves passing through the specimen have limited bandwidth and therefore can be
approximated very well by a low-dimensional linear space spanned by translations of a
well-localized function $u$. Just like in the PRISM algorithm recently published by C. Ophus,
we utilize the linearity of the Schrödinger equation, but perform the approximations with
functions that are well localized in real space instead of Fourier space. This way, we achieve
a similar computational speedup as PRISM, but at a much lower memory consumption and
reduced numerical error due to avoiding virtual copies of the probe waves interfering with the
result. Our approach also facilitates faster recomputations if local changes are made to the specimen
such as changing a single atomic column.
\end{abstract}

\section{Introduction}

The simulation of scanning transmission electron microscopy (STEM) images is frequently
necessary for a correct interpretation of the experimental data.
Unfortunately, STEM image simulation is a computationally highly demanding task. Depending on
specimen complexity and desired resolution, the simulation of a single STEM image
may take many hours or even days to complete \cite{blom18}. If multiple images need to be simulated
for different atomic configurations, this task quickly becomes infeasible.

At its core, the problem is that a partial differential equation (PDE) needs to be solved
many times. For STEM imaging, the relativistically correct PDE describing the transmission of an electron
wave through the electrostatic potential of a specimen is given by the Dirac equation. 
In it, the electron wave is represented by a function with values in a four-dimensional complex space
and its computation is not practical. The commonly used approximation for the electron wave is a single
complex-valued function solving the Schrödinger equation for fast electrons \cite{kirkland10}.

For the purpose of STEM image simulation, the initial conditions to the PDE typically represent 
different positions of the STEM probe, each of which corresponds to a single
pixel in the two dimensional output image. In order to account for thermal diffuse scattering (TDS)
caused by thermal vibrations of the atoms and partial coherence of the STEM probes, the entire simulation
is repeated several times with slightly different probe and specimen parameters and the results are averaged.
In the standard \textit{multislice algorithm} \cite{kirkland10, cowley57} for example,
the Schrödinger equation needs to be solved several times for the computation of one pixel
of the resulting STEM image. Consequently, even for a small output image sampled at $256 \times 256$ pixels, this results in
having to solve the Schrödinger equation millions of
times.

There has been some progress lately on accelerating these computations. A straightforward
approach is to limit the size of the computation window to
a small region around the center of the individual STEM probes. This way, the Schrödinger equation only needs to be solved on
an area that may be much smaller than
the entire specimen area. This provides good approximations to the true
solution due to the relatively fast decay of focused STEM probes and yields a speed up of approximately
$f^2$ compared to the standard multislice algorithm if the size of the simulation window
is reduced by a factor of $f$ in both directions. There is a rumour about MULTEM \cite{lobato15, lobato16} implementing
this approach, but it is not mentioned in the corresponding publications. MULTEM also comprises a
higher order expansion of the multislice
split-step solution and an efficient implementation on the GPU. STEMsalabim \cite{oelerich17}
is a recent implementation of the standard Multislice algorithm that achieves short computation times
via distributed computing using MPI.

The PRISM algorithm \cite{ophus17} and its implementation Prismatic \cite{pryor17, dacosta21} take a
very different approach. Instead of propagating the STEM probes themselves through the specimen,
PRISM utilizes the linearity of the Schrödinger equation and propagates certain exponential functions
through the specimen, which are given by Dirac deltas in Fourier space. The probes are then approximated by
linear combinations of these Dirac deltas, which gives the corresponding
solutions to the Schrödinger equation due to its linearity. One advantage of this approach is the fact
that the frequencies of the STEM probes are limited by the objective aperture, which typically results in a smaller
number of Dirac deltas than STEM probe positions. An additional advantage 
is that it is possible to consider only every $f$-th Dirac delta in $x$ and $y$ directions after the discretization,
thereby reducing by $f^2$ the number of times that the Schrödinger equation needs to be solved at
a negligible loss of accuracy for small $f$.
Since additionally only those exponential functions that correspond to frequencies within the (potentially small)
objective aperture need to be considered,
this method achieves shorter computation times than the simple reduction of the simulation
window size by a factor of $f^2$.

Here, we present a new method based on the idea of PRISM. Instead of using functions that are localized in
Fourier space such as the Dirac deltas of PRISM,
we consider functions that are localized in real space
coordinates for the approximation of the STEM probes.
This has several advantages on which we will elaborate in more detail in this paper: 
(i) it reduces the numerical error in the probe approximation, 
(ii) avoids an excessive memory consumption for large specimen, and 
(iii) reduces the required amount of additional computations if local changes are made to the specimen structure.
At the same time, we are able to keep the benefits of PRISM. In particular, a similar
computational speedup is achieved and we keep the ability to flexibly trade off computation time
against accuracy by means of several parameters described in \cref{rs_approximation}.

There exists a multitude of codes for electron microscopy simulation that serve various purposes, many of which are freely available
online. Aside from the algorithms mentioned above, lists with additional references can be found in
\cite{pryor17, kirkland16, kirkland10}.

\section{Introduction to LMA}
\label{s:2}
Lattice Multislice Algorithm (LMA) utilizes translations of a special initial wave $u$ for which a variant of the \textit{multislice algorithm} (see Subsection \ref{s:msa}) is applied. 
Given a set $\mathcal{I}$ of shift positions $i\in\mathcal{I}$ (considered in Subsection \ref{s:4.2.1}), we obtain a set of initial waves $u_i:=u(\cdot-i)$ that are used to approximate the initial waves of the \textit{probes} $\psi_p^{\mathrm{init}}$ for $p\in\mathcal{P}$. Here $\mathcal{P}$ is the set of probe positions for which we want to calculate  $\hat{\psi}_p^{\mathrm{exit}}$, see Subsection \ref{s:3.1}.
The calculations of the electron wave propagation are performed on a rectangular grid of $X\times Y$ points (see Section \ref{s:4}). 
For convenience, $\mathcal{P}$ can be taken to be a subset of this grid, while $\mathcal{I}$ is often a subset of $\mathcal{P}$. The number of elements of these sets usually satisfy the inequalities  $|\mathcal{I}|<|\mathcal{P}|<XY$, although $\abs{\mathcal{I}} >= \abs{\mathcal{P}}$ is also possible.

In contrast, PRISM approximates the probes via linear combinations of a set $\mathcal{K}$ of elementary waves represented by Dirac deltas in Fourier space (see Subsection \ref{s:4.1}).  
The set $\mathcal{K}$ is uniquely determined by the aperture function $A(k)$ from (\ref{aperture_function}) and the discretization provided by the $X\times Y$ grid in real space. This discretization fixes the scaling of $k$ the Fourier space and determines the number of elementary waves that is limited by the number of discrete points  $k$ satisfying $\|k\|<k_{\max}$.

The linear space $\mathrm{span}(\mathcal{K})$ of the linear combinations of the elementary waves from $\mathcal{K}$ is $|\mathcal{K}|$-dimensional and approximates the probes extremely well. This indicates the feasibility of finding a linear space of dimension $|\mathcal{K}|$ or slightly larger that approximates the probes very well. 
In LMA we represent this space as the linear combinations $\mathrm{span}(u_i, i\in\mathcal{I})$ of $u_i$ using different options for $u$. Standard choices for $u$ are an initial probe wave or a Gaussian, both considered in Subsection \ref{rsapprox_inputwaves}. Other options are also discussed there. One of them establishes a direct theoretical link between PRISM and LMA. It is based on a special representation of $u$ as a tensor product of univariate trigonometric polynomials.

In one dimension, the elementary waves $e^{ikt}$ are the Dirac deltas at the discrete points in Fourier space. The trigonometric polynomials $\mathcal{T}_n$ are the linear subspace of the space of $2\pi$-periodic functions spanned by $e^{ikt}$ for $k=-n,-n+1,...,n$. 
An alternative basis of $\mathcal{T}_n$ is the set of translations $\varphi_n(t-\frac{2\pi j}{2n+1})$ for $j=-n,-n+1,...,n$ of the trigonometric polynomial
\begin{equation}
\label{rauhut}
  \varphi_n(t) = \sum_{k=-n}^n e^{ikt}\cos\left(\frac{k\pi}{2n+2}\right)
\end{equation}
defined in \cite{rauhut05}.
Let $W$ be the $(2n+1)\times(2n+1)$ modulation matrix with entries $W_{k,j}=e^{-\frac{2\pi jk}{2n+1}}$. Then the $(2n+1)\times(2n+1)$ transformation matrix $M$ from the basis of the elementary waves to the translations of the localized trigonometric polynomial from (\ref{rauhut}) is $$M=W\mathrm{diag}\left(\cos\Big(\frac{k\pi}{2n+2}\Big)\right)\,,$$
where $\mathrm{diag}(a_k)$ is the diagonal matrix with entries $a_k$ on the main diagonal. Using that $\frac{1}{2n+1}W$ is a Hermitian matrix, it is easy to see that
$$ M^{-1}=\frac{1}{(2n+1)^2}\mathrm{diag}\left(1\Big/\cos\Big(\frac{k\pi}{2n+2}\Big)\right)W^*\,,
$$
where $W^*=\bar{W}^{\mathrm{T}}$ is the conjugate transpose of $W$.
Matrices $M$ and $M^{-1}$ provide easy transformations from one basis to another in one dimension.

To use the above framework in two dimensions we define $u(x, y) := \varphi_n(x)\varphi_n(y)$ to be the tensor product of localized trigonometric polynomials in each direction. This allows an easy transition from the global elementary waves to localized waves using the matrices $M$ and $M^{-1}$. 
A minor inconvenience is the fact that instead of the set $\mathcal{K}$ of elementary waves related to the points  $k=(k_x,k_y)$ satisfying $\|k\|=\sqrt{k_x^2+k_y^2}\leq n$ we now consider a larger set $\mathcal{K}_+$ related to $k$ satisfying $|k_x|\leq n$ and $|k_y|\leq n$ to receive the set $\mathcal{I}$ of translations $i$ of localized waves $u_i$ with $|\mathcal{I}|=|\mathcal{K}_+|\approx 4n^2$ compared to $|\mathcal{K}|\approx \pi n^2$. 
However, this setup provides more symmetries and less complicated implementation. In particular, the representations of the probes via linear combinations of the translations $u_i$ of $u$ are received from the ones for the elementary waves by multiplications by $M^{-1}$. 
It should be noted that starting from $\mathrm{span}(\mathcal{K}_f)$ one can receive in a similar way a corresponding generating function $u_f$ for a reduced set of translations $\mathcal{I}_f$.

\section{Background}

We briefly summarize the most relevant steps of STEM image formation and the standard Multislice
algorithm and introduce the notation used in the rest of the paper.
A much more detailed description of the standard algorithms for STEM image simulation 
can be found in the book by Kirkland \cite{kirkland10}.

\subsection{STEM image formation}\label{STEM_image_formation}
\label{s:3.1}
In STEM, an image is formed by scanning a specimen one point at a time. The scanning positions
are arranged on a standard lattice and at each of these positions a focused STEM probe is
transmitted through the specimen and recorded subsequently after exiting the specimen. The detector has
the shape of an annulus or a disk and records the squared magnitude of the electron wave
in Fourier space. It is possible to record the squared magnitude of the entire electron wave as
a pixelated image (4D STEM), the integrated values on many small annuli (3D STEM) or the
integrated value over one large annulus or disk (2D STEM).

For a given probe position $p\in\R^2$, the STEM output is generated in three steps as follows.

\paragraph{1.} A focused STEM probe is formed at the position $p$ on the surface of the specimen. The probe wave function
is given by
\begin{equation}\label{probe_wave_function}
  \psi_p^\text{init}(x) = \mathcal{F}^{-1}\big(\exp(-\ii\chi)A\big)(x-p) \qquad\forall\,x\in\R^2,
\end{equation}
where $\mathcal{F}^{-1}$ is the inverse Fourier transform and $\chi:\R^2\rightarrow\R$ is the wave aberration function
\begin{equation}\label{aberration_function}
  \chi(k) = \frac{1}{2}\pi C_s\lambda^3 \norm{k}_2^4 - \pi Z \lambda \norm{k}_2^2 \qquad\forall\,k\in\R^2
\end{equation}
with third order spherical aberration $C_s\in\R$, focus $Z\in\R$ and electron wavelength $\lambda>0$.
The aperture function $A$ acts as an ideal lowpass filter on the STEM probe and is defined as
\begin{equation}\label{aperture_function}
  A(k) = \mathbbm{1}_{B_{\alpha_\text{max}/\lambda}(0)}(k) = \begin{cases} 1, & \lambda\norm{k}_2 < \alpha_\text{max}, \\ 0, & \text{otherwise} \end{cases} \qquad\forall\,k\in\R^2,
\end{equation}
where $\alpha_\text{max}>0$ is the maximum semiangle allowed by the objective aperture. For simplicity,
no higher order or anisotropic aberrations have been included in the wave aberration function. The
precise form of the aberration function used does not affect the results presented in this paper.

\paragraph{2.} The probe is transmitted through the specimen. Although not relativistically
correct, the transmission is commonly modeled by the Schrödinger equation for fast electrons,
\begin{equation}\label{schroedinger_equation}
  \frac{\partial\psi(r)}{\partial z} = \frac{\ii\lambda}{4\pi}\Delta_{xy}\psi(r) + \ii\sigma V(r)\psi(r) \qquad\forall\,r=(x,y,z)\in\R^3,
\end{equation}
with $\psi_p^\text{init}$ from step 1 as the initial condition at the entrance plane of the specimen,
\begin{equation*}
  \psi(x,y,0) = \psi_p^\text{init}(x,y) \qquad\forall\,(x,y)\in\R^2.
\end{equation*}
Here, $\Delta_{xy} = \frac{\partial^2}{\partial x^2} + \frac{\partial^2}{\partial y^2}$
is the two dimensional Laplace operator, $\sigma>0$ is the beam sample interaction constant and
$V:\R^3\rightarrow\R_{>0}$ is the electrostatic potential of the specimen.
The electron wave at the exit plane of the specimen, also called exit wave, is denoted by $\psi_p^\text{exit}(x,y) := \psi_p(x,y,t)$,
where $\psi_p$ is the solution to the Schrödinger equation with initial condition $\psi_p^\text{init}$ and
$t>0$ is the specimen thickness.

\paragraph{3.} A detector records the squared amplitude of the exit wave in Fourier space.
The final STEM output for the probe position $p$ is then given as
\begin{align*}
  &\text{4D STEM: }\quad \left(\Big|\mathcal{F}\big(\psi_p^\text{exit}\big)(ad_x+bd_y)\Big|^2\right)_{\substack{a=-A,\ldots,A \\ b=-B,\ldots,B}}  \tag*{$(A,B\in\N,\;d_x,d_y>0)$} \\
  &\text{3D STEM: }\quad \left(\int_{B_{(a+1)r}(0)\backslash B_{ar}(0)} \abs{\mathcal{F}\big(\psi_p^\text{exit}\big)(k)}^2\dk\right)_{a=0,\ldots,A} \tag*{$(A\in\N,\;r>0)$} \\
  &\text{2D STEM: }\quad \int_{B_{r_2}(0)\backslash B_{r_1}(0)} \abs{\mathcal{F}\big(\psi_p^\text{exit}\big)(k)}^2\dk, \tag*{$(r_1, r_2\ge 0)$}
\end{align*}
where $\mathcal{F}$ is the (forward) Fourier transform.\\

The collection of the pixel values in step 3 for all probe positions $p$
forms the 2D, 3D or 4D STEM image. The probe positions generally lie on a rectangular standard
lattice and are contained within a rectangular domain corresponding to the specimen area
of which a STEM image is simulated, referred to as the simulation domain in the following.
The simulation domain is usually assumed to be periodic because of the use of the discrete Fourier
transform (DFT) in the Multislice algorithm.

The description of STEM image formation given above assumes a perfectly coherent illumination and a
completely static specimen. In reality, the illumination in a STEM is only partially
coherent and the atoms in the specimen vibrate slightly so that the
specimen potential $V$ is not constant for the entire duration of image acquisition. In order to
account for these effects, steps 1 to 3 are typically repeated several times with slightly
different probe wave functions $\psi_p^\text{init}$ and specimen potentials $V$. In the end,
these results are averaged to compute the final STEM image.

\subsection{The Multislice Algorithm}
\label{s:msa}
Steps 1 and 3 above can be easily computed, while solving the Schrödinger equation in
step 2 is computationally more demanding. The Schrödinger equation is generally solved with
the Multislice algorithm, which is a split-step method similar to the Euler method.
Another popular method, the Bloch wave method, requires solving a large eigenvalue/eigenvector
problem and is feasible only for very small specimen \cite{kirkland10}.

\begin{figure}[h]
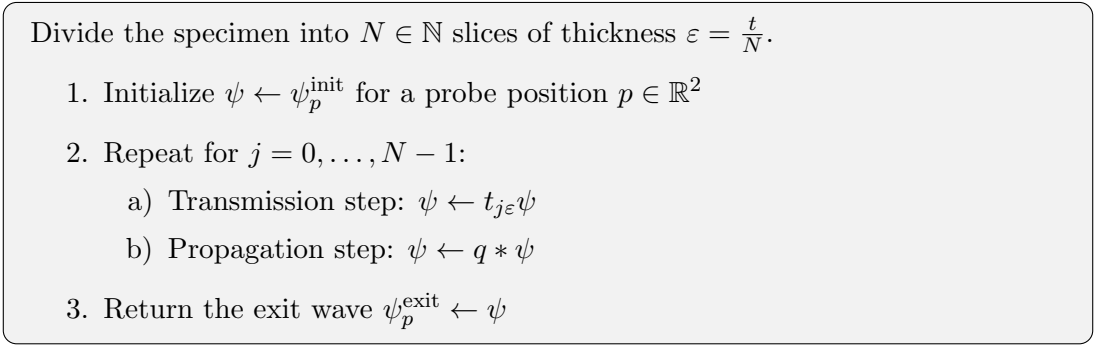

  \hfill
  \begin{mdframed}[backgroundcolor=BoxBg, skipabove=10pt, skipbelow=10pt, leftmargin=5pt, rightmargin=5pt, roundcorner=5pt]
  \begin{minipage}{0.98\textwidth}
    Divide the specimen into $N\in\N$ slices of thickness $\varepsilon = \frac{t}{N}$.
    \begin{enumerate}
      \item Initialize $\psi \leftarrow \psi_p^\text{init}$ for a probe position $p\in\R^2$
      \item Repeat for $j=0,\ldots,N-1$:
            \begin{enumerate}
              \item Transmission step: $\psi \leftarrow t_{j\varepsilon}\psi$
              \item Propagation step: $\psi \leftarrow q \ast \psi$
            \end{enumerate}
      \item Return the exit wave $\psi_p^\text{exit} \leftarrow \psi$
    \end{enumerate}
  \end{minipage}
  \end{mdframed}
  \caption{Pseudocode for solving the Schrödinger equation with the Multislice algorithm.
           The error made for each slice is in $\mathcal{O}(\varepsilon^2)$ so that the
           overall error is in $\mathcal{O}(\varepsilon)$ if $N \sim \frac{1}{\varepsilon}$.}\label{multislice_pseudocode}
\end{figure}

In the Multislice algorithm, the specimen is divided into thin slices along the beam
direction $z$ and the Schrödinger equation \bref{schroedinger_equation} is split into two parts,
the propagation term $\frac{\ii\lambda}{4\pi}\Delta_{xy}\psi(r)$ and the transmission term
$\ii\sigma V(r)\psi(r)$. This yields an approximate solution to the Schrödinger equation by
alternately multiplying the electron wave with the transmission function
\begin{equation*}
  t_z(x,y) := \exp\bigg(\ii\sigma\underbrace{\int_z^{z+\varepsilon} V(x,y,z')\dz'}_{=:v_\varepsilon(x,y,z)}\bigg) = \exp\big(\ii\sigma v_\varepsilon(x,y,z)\big) \qquad\forall\,x,y\in\R^2
\end{equation*}
and convolving it with the propagation function
\begin{equation*}
  q(x,y) := \frac{-\ii}{\lambda\varepsilon}\exp\left(\frac{\ii\pi}{\lambda\varepsilon}(x^2+y^2)\right) \qquad\forall\,x,y\in\R^2,
\end{equation*}
where $\varepsilon>0$ is the slice thickness. The Multislice algorithm is summarized in
\cref{multislice_pseudocode}.

\paragraph{Complexity:} In order to estimate the complexity of the Multislice algorithm, we assume that the electron waves
are discretized as complex images on a grid of $X\times Y$ pixels. The propagation step may
either be computed as a convolution in real space or as a pointwise multiplication in
Fourier space, using the convolution theorem $q\ast\psi = \mathcal{F}^{-1}\big(\mathcal{F}(q)\mathcal{F}(\psi)\big)$ and
the fast Fourier transform. The Fourier transform $\mathcal{F}(q)$ is
a constant and can be precomputed. Performing the Multislice algorithm for a single probe position
$p$ requires approximately
\begin{equation*}
  T_{\text{Multislice, XY}} = N \big(2XY + 2XY\log(XY)\big)
\end{equation*}
operations if the propagation step is performed in Fourier space and
\begin{equation*}
  \tilde T_{\text{Multislice, XY}} = N (XY + XYK_1K_2)
\end{equation*}
operations if the propagation step is performed as a convolution in real space, where the
propagation function $q$ is discretized to $K_1\times K_2$ pixels.

We remark that it is not always necessary to perform the Multislice computations on the
entire simulation domain. Both the Fourier space and the real space variant can be computed
on smaller windows adjusted to the extent of the probe wave function. However, caution must
be exercised when using the Fourier space variant of the Multislice algorithm because of the
periodicity of the DFT. If $X$ and $Y$ are chosen sufficiently large, this problem is alleviated
by the decay of the STEM probe to zero. A too small value of $X$ and $Y$ on the other hand
may cause a significant approximation error due to the wrap-around effect of the DFT.

\section{Approximation of STEM probes}
\label{s:4}
Solving the Schrödinger equation is the computational bottleneck of STEM image simulation.
From a purely mathematical point of view, there are two complementary ways to alleviate
this problem. On the one hand, one can try to find better ways of solving the Schrödinger
equation, which means using a faster or more accurate method than the Multislice algorithm.
On the other hand, one can try to reduce the number of times that the Schrödinger
equation needs to be solved. In this section, we will be discussing the latter approach.

In the following, the size of the discretized simulation domain in pixels is given by
$X\times Y$ for $X, Y\in\N$.

\subsection{Fourier space approximation in PRISM}
\label{s:4.1}
The PRISM algorithm \cite{ophus17} is based on one key observation: instead of propagating
the STEM probes themselves through the specimen, it is equally possible to first propagate
each of their discrete frequencies individually and afterwards reconstruct the STEM probes
as linear combinations of their frequencies.

Consider the discretized STEM probes $\hat\psi_p^\text{init}$ on a grid
of $X\times Y$ pixels with a total physical dimension of $l_x \times l_y$ {\AA}ngströms, where $l_x,l_y\in\R$. According to \cref{probe_wave_function}, the
Fourier space values $\mathcal{F}\big(\hat\psi_p^\text{init}\big)\in\C^{X\times Y}$ may be given as
\begin{equation*}
  \mathcal{F}\left(\hat\psi_p^\text{init}\right)_{k_x, k_y} = \exp\big(-\ii\chi\big(g(k_x,k_y)\big)\big)A\big(g(k_x,k_y)\big)\mu_{-p}\big(g(k_x,k_y)\big)
\end{equation*}
for all $k_x\in\{1,\ldots,X\}$, $k_y\in\{1,\ldots,Y\}$, where
\begin{equation*}
  g(k_x,k_y) = \left(\frac{k_x-1 - \floor{X/2}}{l_x}, \frac{k_y-1 - \floor{Y/2}}{l_y}\right)
\end{equation*}
converts pixel coordinates $(k_x, k_y)$ to the corresponding frequencies. The Fourier transform turns translations in
real space coordinates into modulations in Fourier space, where the modulation by $p\in\R^2$
is defined as
\begin{equation*}
  \mu_p: \R^2 \rightarrow \C, \quad k \mapsto \exp(2\pi\ii k\cdot p).
\end{equation*}
The value $\mathcal{F}\left(\hat\psi_p^\text{init}\right)_{k_x, k_y}$ is nonzero if and only if
\begin{equation*}
  (k_x, k_y) \in \mathcal{K} := \big\{ (k_x, k_y) \in \{1,\ldots, X\} \times \{1, \ldots, Y\} \mid A\big(g(k_x, k_y)\big) \neq 0 \big\}.
\end{equation*}

Expressing $\mathcal{F}\big(\hat\psi_p^\text{init}\big)$ as a linear combination of discrete dirac deltas, it follows that
\begin{equation*}
  \mathcal{F}\left(\hat\psi_p^\text{init}\right) = \sum_{(k_x', k_y')\in\mathcal{K}} \mathcal{F}\left(\hat\psi_p^\text{init}\right)_{k_x', k_y'} \delta_{(k_x', k_y')},
\end{equation*}
where $\delta_b\in\R^{X\times Y}$ with
\begin{equation*}
  (\delta_b)_a = \begin{cases} 1, & \text{if } a = b, \\ 0, & \text{otherwise} \end{cases}
\end{equation*}
for all $a,b\in\{1,\ldots,X\} \times \{1,\ldots,Y\}$. Consequently, taking the inverse discrete Fourier transform yields
\begin{equation}\label{prism_probe_frequencies}
  \hat\psi_p^\text{init} = \sum_{(k_x', k_y')\in\mathcal{K}} \mathcal{F}\left(\hat\psi_p^\text{init}\right)_{k_x', k_y'} \underbrace{\mathcal{F}^{-1}\left(\delta_{(k_x', k_y')}\right)}_{=:w_{k_x', k_y'}}.
\end{equation}
The plane waves $w_{k_x', k_y'}$ are the functions that are transmitted through the
specimen in the PRISM algorithm instead of the probes $\hat\psi_p^\text{init}$.

Let $\mathcal{P}\subseteq\R^2$ be the set of all probe positions that form the final STEM image.
If $W_{k_x', k_y'}$ denotes the solution to the Schrödinger equation for the initial condition
$w_{k_x', k_y'}$, then the exit waves are given by
\begin{equation}\label{prism_linearcombination_exit}
  \hat\psi_p^\text{exit} = \sum_{(k_x', k_y')\in\mathcal{K}} \mathcal{F}\left(\hat\psi_p^\text{init}\right)_{k_x', k_y'} W_{k_x', k_y'}
\end{equation}
for all $p\in \mathcal{P}$ due to the linearity of the Schrödinger equation.

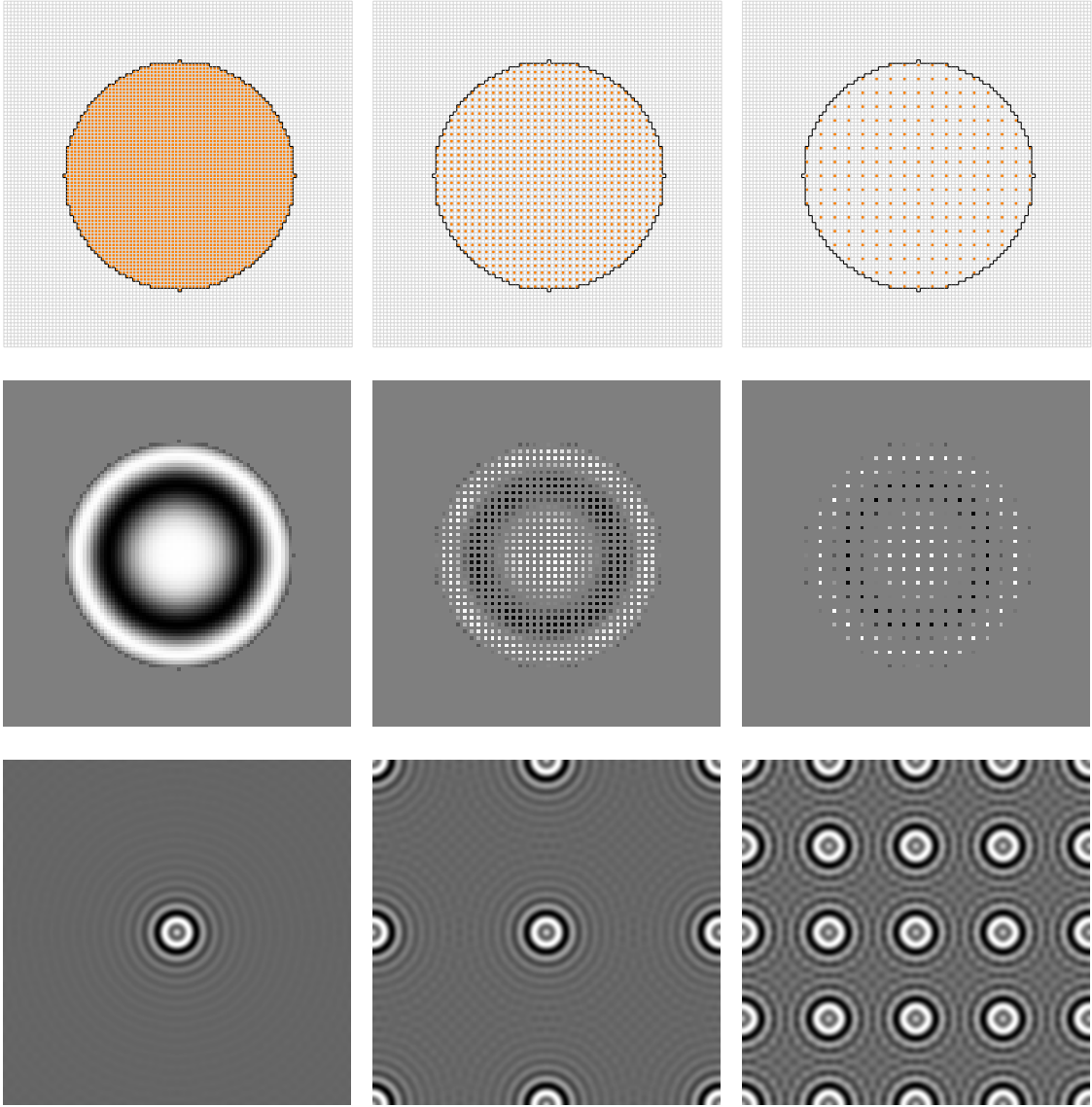
\begin{figure}[!t]
  \begin{center}
  \begin{minipage}{0.32\textwidth}
  \begin{center}
    \begin{tikzpicture}[yscale=-1, x=0.1225em, y=0.1225em]
      \draw[gray!30](0, 0) grid[step=1] (100, 100);
      \foreach \X in {0, 1, ..., 99} {
        \foreach \Y in {0, 1, ..., 99} {
          \pgfmathtruncatemacro{\itest}{(\X-50)*(\X-50) + (\Y-50)*(\Y-50)}
          \ifnum\itest<1090
            \draw[gray!30,fill=orange] (\X, \Y) rectangle (\X+1, \Y+1);
            
            \pgfmathtruncatemacro{\itestx}{(\X+1-50)*(\X+1-50) + (\Y-50)*(\Y-50)}
            \ifnum\itestx>1089
              \draw[black] (\X+1, \Y) -- (\X+1, \Y+1);
            \fi
            \pgfmathtruncatemacro{\itestx}{(\X-1-50)*(\X-1-50) + (\Y-50)*(\Y-50)}
            \ifnum\itestx>1089
              \draw[black] (\X, \Y) -- (\X, \Y+1);
            \fi
            \pgfmathtruncatemacro{\itestx}{(\X-50)*(\X-50) + (\Y+1-50)*(\Y+1-50)}
            \ifnum\itestx>1089
              \draw[black] (\X, \Y+1) -- (\X+1, \Y+1);
            \fi
            \pgfmathtruncatemacro{\itestx}{(\X-50)*(\X-50) + (\Y-1-50)*(\Y-1-50)}
            \ifnum\itestx>1089
              \draw[black] (\X, \Y) -- (\X+1, \Y);
            \fi
            
          \fi
        }
      }
    \end{tikzpicture}
  \end{center}
  \end{minipage}
  \hfill
  \begin{minipage}{0.32\textwidth}
  \begin{center}
    \begin{tikzpicture}[yscale=-1, x=0.1225em, y=0.1225em]
      \draw[gray!30](0, 0) grid[step=1] (100, 100);
      \foreach \X in {0, 1, ..., 99} {
        \foreach \Y in {0, 1, ..., 99} {
          \pgfmathtruncatemacro{\itest}{(\X-50)*(\X-50) + (\Y-50)*(\Y-50)}
          \ifnum\itest<1090
            \pgfmathparse{mod(\X, 2)} \pgfmathtruncatemacro{\htest}{\pgfmathresult}
            \pgfmathparse{mod(\Y, 2)} \pgfmathtruncatemacro{\vtest}{\pgfmathresult}
            \ifnum\htest=0
              \ifnum\vtest=0
                \draw[gray!30,fill=orange] (\X, \Y) rectangle (\X+1, \Y+1);
              \fi
            \fi
            
            \pgfmathtruncatemacro{\itestx}{(\X+1-50)*(\X+1-50) + (\Y-50)*(\Y-50)}
            \ifnum\itestx>1089
              \draw[black] (\X+1, \Y) -- (\X+1, \Y+1);
            \fi
            \pgfmathtruncatemacro{\itestx}{(\X-1-50)*(\X-1-50) + (\Y-50)*(\Y-50)}
            \ifnum\itestx>1089
              \draw[black] (\X, \Y) -- (\X, \Y+1);
            \fi
            \pgfmathtruncatemacro{\itestx}{(\X-50)*(\X-50) + (\Y+1-50)*(\Y+1-50)}
            \ifnum\itestx>1089
              \draw[black] (\X, \Y+1) -- (\X+1, \Y+1);
            \fi
            \pgfmathtruncatemacro{\itestx}{(\X-50)*(\X-50) + (\Y-1-50)*(\Y-1-50)}
            \ifnum\itestx>1089
              \draw[black] (\X, \Y) -- (\X+1, \Y);
            \fi
            
          \fi
        }
      }
    \end{tikzpicture}
  \end{center}
  \end{minipage}
  \hfill
  \begin{minipage}{0.32\textwidth}
  \begin{center}
    \begin{tikzpicture}[yscale=-1, x=0.1225em, y=0.1225em]
      \draw[gray!30](0, 0) grid[step=1] (100, 100);
      \foreach \X in {0, 1, ..., 99} {
        \foreach \Y in {0, 1, ..., 99} {
          \pgfmathtruncatemacro{\itest}{(\X-50)*(\X-50) + (\Y-50)*(\Y-50)}
          \ifnum\itest<1090
            \pgfmathparse{mod(\X+2, 4)} \pgfmathtruncatemacro{\htest}{\pgfmathresult}
            \pgfmathparse{mod(\Y+2, 4)} \pgfmathtruncatemacro{\vtest}{\pgfmathresult}
            \ifnum\htest=0
              \ifnum\vtest=0
                \draw[gray!30,fill=orange] (\X, \Y) rectangle (\X+1, \Y+1);
              \fi
            \fi
            
            \pgfmathtruncatemacro{\itestx}{(\X+1-50)*(\X+1-50) + (\Y-50)*(\Y-50)}
            \ifnum\itestx>1089
              \draw[black] (\X+1, \Y) -- (\X+1, \Y+1);
            \fi
            \pgfmathtruncatemacro{\itestx}{(\X-1-50)*(\X-1-50) + (\Y-50)*(\Y-50)}
            \ifnum\itestx>1089
              \draw[black] (\X, \Y) -- (\X, \Y+1);
            \fi
            \pgfmathtruncatemacro{\itestx}{(\X-50)*(\X-50) + (\Y+1-50)*(\Y+1-50)}
            \ifnum\itestx>1089
              \draw[black] (\X, \Y+1) -- (\X+1, \Y+1);
            \fi
            \pgfmathtruncatemacro{\itestx}{(\X-50)*(\X-50) + (\Y-1-50)*(\Y-1-50)}
            \ifnum\itestx>1089
              \draw[black] (\X, \Y) -- (\X+1, \Y);
            \fi
            
          \fi
        }
      }
    \end{tikzpicture}
  \end{center}
  \end{minipage}
  \end{center}
  
  \begin{center}
  \includegraphics[width=0.32\textwidth, keepaspectratio]{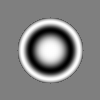} \hfill
  \includegraphics[width=0.32\textwidth, keepaspectratio]{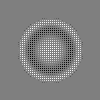} \hfill
  \includegraphics[width=0.32\textwidth, keepaspectratio]{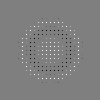}
  \end{center}
  
  \begin{center}
  \includegraphics[width=0.32\textwidth, keepaspectratio]{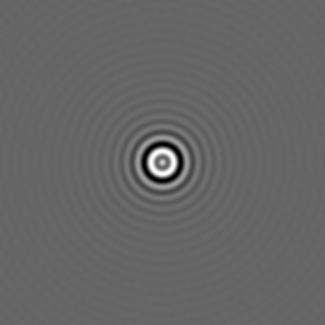} \hfill
  \includegraphics[width=0.32\textwidth, keepaspectratio]{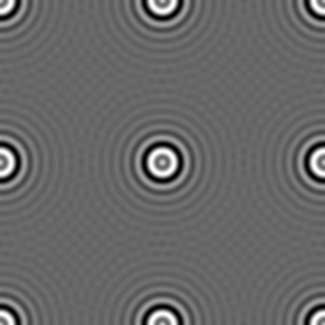} \hfill
  \includegraphics[width=0.32\textwidth, keepaspectratio]{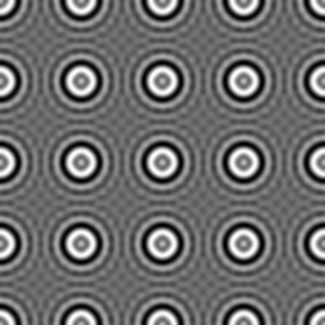}
  \end{center}
  \caption{From left to right: the approximation \bref{prism_approximation_init} with
           PRISM interpolation factor $f=1$, $f=2$ and $f=4$. From top to bottom: the
           elements of the set $\mathcal{K}_f$ denoted by the filled pixels (top row),
           the real part of \bref{prism_approximation_init} in Fourier space, cropped to a
           small area around the non-zero frequencies, (middle row)
           and the real part of \bref{prism_approximation_init} in real space (bottom row). Note that the images have
           periodic boundary conditions due to the periodicity of the discrete Fourier
           transform.}\label{prism_approximation_figure}
\end{figure}

Although the number of plane waves $w_{k_x', k_y'}$ may be much smaller than the number of
probe positions $\abs{\mathcal{P}}$ and the Schrödinger equation only needs to be solved for the plane wave initial
conditions, this approach may still be slower than the standard Multislice algorithm for
two reasons:
\begin{enumerate}
  \item Evaluating \bref{prism_linearcombination_exit} for every probe position $p\in \mathcal{P}$ requires a
        large number of floating point operations and may take a long time depending on the
        size of $\abs{\mathcal{K}}$.
  \item All of the $W_{k_x', k_y'}$ for $(k_x', k_y')\in\mathcal{K}$ need to be kept in memory
        at the same time, causing large memory consumption and potentially slow memory access.
\end{enumerate}
To a certain extent, both of these problems are alleviated by the PRISM interpolation factor.
This is a factor $f\in\N$ which reduces the number of plane waves by considering only
every $f$-th dirac delta in Fourier space. Formally, the set $\mathcal{K}$ is replaced by
\begin{equation*}
  \mathcal{K}_f := \mathcal{K} \cap \left(\left(f\Z + 1 + \left\lfloor\frac{X}{2}\right\rfloor\right) \times \left(f\Z + 1 + \left\lfloor\frac{Y}{2}\right\rfloor\right)\right)
\end{equation*}
and the approximation
\begin{equation}\label{prism_approximation_init}
  \hat\psi_p^\text{init} \approx \sum_{(k_x', k_y')\in\mathcal{K}_f} \mathcal{F}\left(\hat\psi_p^\text{init}\right)_{k_x', k_y'} w_{k_x', k_y'}
\end{equation}
leads to
\begin{equation}\label{prism_approximation_exit}
  \hat\psi_p^\text{exit} \approx \sum_{(k_x', k_y')\in\mathcal{K}_f} \mathcal{F}\left(\hat\psi_p^\text{init}\right)_{k_x', k_y'} W_{k_x', k_y'}.
\end{equation}
The set $\mathcal{K}_f$ is depicted in the top row of \cref{prism_approximation_figure} for $f\in\{1,2,4\}$.
This is an elegant way to speed up the computations and reduces the number of times that the
Schrödinger equation needs to be solved by a factor of $\abs{\mathcal{K}}/\abs{\mathcal{K}_f} \approx f^2$.
The PRISM algorithm is summarized in \cref{prism_pseudocode}.

\begin{figure}[h]
  \hfill
  \begin{mdframed}[backgroundcolor=BoxBg, skipabove=10pt, skipbelow=10pt, leftmargin=5pt, rightmargin=5pt, roundcorner=5pt]
  \begin{minipage}{0.98\textwidth}
    Choose an interpolation factor $f\in\N$.
    \begin{enumerate}
      \item For each $k'\in\mathcal{K}_f$, use the Multislice algorithm (cf. \cref{multislice_pseudocode})
            to calculate the solution $W_{k'}$ to the Schrödinger equation with initial condition $w_{k'}$.
      \item For each probe position $p\in \mathcal{P}$, calculate an approximation to the corresponding
            exit wave using \cref{prism_approximation_exit}.
    \end{enumerate}
  \end{minipage}
  \end{mdframed}
  \caption{Pseudocode for the PRISM algorithm.}\label{prism_pseudocode}
\end{figure}

A major practical difficulty in this approach is the implicit repetition of the STEM probe
in the approximation \bref{prism_approximation_init}. Using an interpolation factor of $f\in\N$ causes
$f^2$ evenly spaced STEM probes to appear in the approximation \bref{prism_approximation_init}, which
means that the approximate solution to the Schrödinger equation in \bref{prism_approximation_exit}
corresponds to propagating $f^2$ copies of the STEM probe through the specimen simultaneously.
This effect is depicted in the bottom row of \cref{prism_approximation_figure}.

Consequently, the approximations \bref{prism_approximation_init} and \bref{prism_approximation_exit} are
valid only locally near the center $p$ of the probe and "interference" with nearby artificial copies of the
probe may introduce significant numerical errors if $f$ is too large. For this reason, the approximations
\bref{prism_approximation_exit} are cropped to a window of $(X/f) \times (Y/f)$ pixels centered at the
probe position in the PRISM algorithm, which is both beneficial and detrimental. It is
beneficial, because it reduces the number of operations that are necessary to evaluate
\bref{prism_approximation_exit} by $f^2$. At the same time, it is detrimental because the useful
area of the simulation window is unavoidably reduced by the same factor $f^2$.

Overall, the number of necessary floating point operations to
compute one STEM image with the PRISM algorithm can be roughly estimated as
\begin{align*}
  T_\text{PRISM} &= \abs{\mathcal{K}_f} T_{\text{Multislice}, XY} + \abs{\mathcal{P}}\abs{\mathcal{K}_f}\frac{XY}{f^2} \\
  &\approx \frac{1}{f^2}\abs{\mathcal{K}} N (2XY + 2XY\log(XY)) + \frac{1}{f^4}\abs{\mathcal{P}}\abs{\mathcal{K}}XY.
\end{align*}
The first term consists of the number of operations necessary to compute the solutions $W_k$ to the
Schrödinger equation with initial condition $w_k$ for all $k\in\mathcal{K}$ and the second term
accounts for the evaluation of \bref{prism_approximation_exit}. Thermal diffuse scattering and partial
coherence are not included in this formula and would result in multipling $T_\text{PRISM}$ by the number
of frozen phonon passes and probe configurations respectively.


\subsection{Real space approximation}\label{rs_approximation}

Since the implicit probe repetition in the PRISM approximation renders large values of
$f$ unusable, the number of solutions $W_{k_x', k_y'}$ that need to be calculated may still be
too large to keep all of them in computer memory at the same time. Additionally, all of the
plane waves $w_{k_x', k_y'}$ have unbounded support in real space coordinates and
span the entire specimen area. Thus, if small local changes are made to the specimen in question,
all of the solutions $W_{k_x', k_y'}$ need to be recomputed.

In order to address these problems, we first take an algebraic point of view on the approximation
of the STEM probes.
Denote by
\begin{equation*}
  \mathcal{P} := \left\{a\begin{pmatrix*} p_x \\ 0 \end{pmatrix*} + b\begin{pmatrix*} 0 \\ p_y \end{pmatrix*} \;\middle|\; a\in\{0,\ldots,P_x-1\}, \; b\in\{0,\ldots,P_y-1\} \right\}
\end{equation*}
the set of all STEM probe positions that form the final STEM image, where
$P_x, P_y\in\N$ are the number of pixels of the output STEM image and $p_x,p_y\in\R$ is
the spacing in between the pixels. The
vector space generated by all probe functions is
\begin{equation*}
  V_\mathcal{P} := \big\langle \hat\psi_p^\text{init} \mid p \in \mathcal{P} \big\rangle_{\C} \subseteq \C^{X\times Y}.
\end{equation*}
By \bref{prism_probe_frequencies}, this vector space satisfies
\begin{equation*}
  V_\mathcal{P} \subseteq \left\langle w_k \mid k \in \mathcal{K} \right\rangle_{\C}
\end{equation*}
and therefore the dimension of $V_\mathcal{P}$ as a complex vector space is bounded by $\abs{\mathcal{K}}$,
regardless of the number of probe positions $\abs{\mathcal{P}}$. The goal is now to find a new basis for $V_\mathcal{P}$,
which avoids the problems of the plane wave basis $(w_k)_{k\in\mathcal{K}}$. This is summarized as
follows:

\begin{mdframed}[backgroundcolor=BoxBg, skipabove=10pt, skipbelow=10pt, leftmargin=5pt, rightmargin=5pt, roundcorner=5pt]
  \vspace*{-1em}
  \begin{task}\label{task_realspace_basis}
    Find a basis $(u_i)_{i\in\mathcal{I}}$ for a vector space $U\subseteq\C^{X\times Y}$ with $V_\mathcal{P}\subseteq U$
    and dimension $\dim_{\C}(U) \le \abs{\mathcal{K}}$ such that the basis elements $u_i$ have a small
    support in real space coordinates.
  \end{task}
\end{mdframed}

Given such a basis $(u_i)_{i\in\mathcal{I}}$, one can first compute solutions to the Schrödinger
equation with the initial conditions $u_i$ and then express the STEM probes $\hat\psi_p^\text{init}$
as linear combinations of the basis elements $(u_i)_{i\in\mathcal{I}}$. Due to the linearity of the Schrödinger
equation, the exit waves $\hat\psi_p^\text{exit}$ can subsequently be computed similarly
as in PRISM.

The requirement on the basis elements $u_i$ to have small support in real space is essential and can
equivalently be phrased by requiring the $u_i$ to be localized in real space. Because of this, a position in
real space can be assigned to each $u_i$ in a similar fashion as is naturally done for the STEM probes.

\subsubsection{The Lattice Multislice Algorithm}\label{lma_description}
\label{s:4.2.1}
We begin by giving a high-level description of a new algorithm for STEM image simulation, which is based on the approximation
of STEM probes by elements $u_i$ as in \cref{task_realspace_basis}. The details are then discussed in the
subsequent sections. It is pointed out that the derivation of this algorithm involves
some simplifications to facilitate its practical implementation. These simplifications
seem natural, but there are likely other approaches that perform equally well of better,
although it is unclear to what extent.\newline

We call the rectangular
standard lattice that is composed of the probe positions $\mathcal{P}$ the probe lattice.
Due to the regular spacing of the STEM probes, it seems reasonable to choose the functions $u_i$
to lie on a rectangular standard lattice too, which we call the input wave lattice. The functions
$u_i$ are called input waves in the following.
We only consider the case where all $u_i$ are equal to a shifted version of one common base
function $u$, just like all STEM probes $\hat\psi_p^\text{init}$ are equal up to their position
$p\in \mathcal{P}$. Without loss of generality, let $\mathcal{I}$ be the set of input wave positions such
that $u_i = u \circ t_i$ for all $i\in\mathcal{I}$ with the translation $t_i: \R^2 \rightarrow \R^2,\; x \mapsto x + i$.
Then
\begin{equation*}
  \mathcal{I} = \left\{a\begin{pmatrix*} r_x \\ 0 \end{pmatrix*} + b\begin{pmatrix*} 0 \\ r_y \end{pmatrix*} + \mathcal{I}_\text{offset} \;\middle|\; a\in\{0,\ldots,R_x-1\}, b\in\{0,\ldots,R_y-1\} \right\}
\end{equation*}
with the number of input waves $R_x, R_y\in\N$ for each direction and
the spacing $r_x, r_y\in\R$ as well as an offset $\mathcal{I}_\text{offset}\in\R^2$, which may also be zero.

Both the probe lattice and the input wave lattice are required to cover the entire simulation domain,
respecting its periodicity. If the simulation domain has a physical size of $l_x \times l_y$ with $l_x, l_y\in\R$,
then this means that
\begin{align}\label{lma_lattices_cover_specimen}
\begin{split}
  &P_xp_x - l_x = P_yp_y - l_y = 0, \\
  &R_xr_x - l_x - (\mathcal{I}_\text{offset})_1 = R_yr_y - l_y - (\mathcal{I}_\text{offset})_2 = 0.
\end{split}
\end{align}

Finally, we also require the probe and input wave lattices to be compatible in the sense that
there exist integers $c,d\ge 1$ such that
\begin{equation}\label{lma_approximation_case1}\tag{Case 1}
  R_x = cP_x \quad\text{and}\quad R_y = dP_y
\end{equation}
or
\begin{equation}\label{lma_approximation_case2}\tag{Case 2}
  P_x = cR_x \quad\text{and}\quad P_y = dR_y.
\end{equation}
In the first case, the number of input waves is greater than or equal to the number of probe positions,
$\abs{\mathcal{I}} \ge \abs{\mathcal{P}}$, and in the second case this inequality is reversed. The key
restriction here is that the dimensions of one lattice are an integral multiple of the dimensions of
the other lattice. While in PRISM the coefficients for the best $L^2$-approximation of the STEM probe
as linear combinations of the $w_k$ for $k\in\mathcal{K}_f$ are given directly by the Fourier space values
of the STEM probe (cf. equations [\ref{prism_probe_frequencies}],[\ref{prism_approximation_init}]), it
is not as straightforward to find coefficients to approximate $\hat\psi_p^\text{init}$ by a linear combination
of the $u_i$. In practice, these coefficients can be precomputed with a linear least squares
approximation; requiring the lattices to be compatible greatly reduces the computational cost of this
precomputation step. This is discussed in more detail in \cref{rsapprox_coefficients}.

\begin{figure}[h]
  \begin{center}
    \begin{tikzpicture}[yscale=-1, x=2.5em, y=2.5em]
      \draw [draw=black] (0, 0) rectangle (9, 7);
      \draw (0, 0) circle (2pt); \draw (1, 0) circle (2pt); \draw (2, 0) circle (2pt); \draw (3, 0) circle (2pt); \draw (4, 0) circle (2pt); \draw (5, 0) circle (2pt); \draw (8, 0) circle (2pt);
      \draw (0, 1) circle (2pt); \draw (1, 1) circle (2pt); \draw (2, 1) circle (2pt); \draw (3, 1) circle (2pt); \draw (4, 1) circle (2pt); \draw (5, 1) circle (2pt); \draw (8, 1) circle (2pt);
      \draw (0, 2) circle (2pt); \draw (1, 2) circle (2pt); \draw (2, 2) circle (2pt); \draw (3, 2) circle (2pt); \draw (4, 2) circle (2pt); \draw (5, 2) circle (2pt); \draw (8, 2) circle (2pt);
      \draw (0, 3) circle (2pt); \draw (1, 3) circle (2pt); \draw (2, 3) circle (2pt); \draw (3, 3) circle (2pt); \draw (4, 3) circle (2pt); \draw (5, 3) circle (2pt); \draw (8, 3) circle (2pt);
      \draw (0, 6) circle (2pt); \draw (1, 6) circle (2pt); \draw (2, 6) circle (2pt); \draw (3, 6) circle (2pt); \draw (4, 6) circle (2pt); \draw (5, 6) circle (2pt); \draw (8, 6) circle (2pt);
      
      \draw (5.5, 0) node[anchor=west] {$\ldots$};
      \draw (5.5, 1) node[anchor=west] {$\ldots$};
      \draw (5.5, 2) node[anchor=west] {$\ldots$};
      \draw (5.5, 3) node[anchor=west] {$\ldots$};
      \draw (5.5, 6) node[anchor=west] {$\ldots$};
      
      \draw (0, 3) node[anchor=north] {$\vdots$};
      \draw (1, 3) node[anchor=north] {$\vdots$};
      \draw (2, 3) node[anchor=north] {$\vdots$};
      \draw (3, 3) node[anchor=north] {$\vdots$};
      \draw (4, 3) node[anchor=north] {$\vdots$};
      \draw (5, 3) node[anchor=north] {$\vdots$};
      \draw (8, 3) node[anchor=north] {$\vdots$};
      
      \draw (5.5, 3) node[anchor=north west] {$\ddots$};
      
      \draw (0.5, 0.5) node[cross] {}; \draw (2.5, 0.5) node[cross] {}; \draw (4.5, 0.5) node[cross] {}; \draw (7.5, 0.5) node[cross] {};
      \draw (0.5, 2.5) node[cross] {}; \draw (2.5, 2.5) node[cross] {}; \draw (4.5, 2.5) node[cross] {}; \draw (7.5, 2.5) node[cross] {};
      \draw (0.5, 5.5) node[cross] {}; \draw (2.5, 5.5) node[cross] {}; \draw (4.5, 5.5) node[cross] {}; \draw (7.5, 5.5) node[cross] {};
      
      \draw [<->] (0, 7.5) -- (9, 7.5);
      \draw (4.5, 7.5) node[anchor=north] {$l_x$};
      
      \draw [<->] (9.5, 0) -- (9.5, 7);
      \draw (9.5, 3.5) node[anchor=west] {$l_y$};
      
      \draw [<->] (0, -0.5) -- (1, -0.5);
      \draw (0.5, -0.5) node[anchor=south] {$p_x$};
      
      \draw [<->] (-0.5, 0) -- (-0.5, 1);
      \draw (-0.5, 0.5) node[anchor=east] {$p_y$};
      
      \draw [<->] (2.5, -0.5) -- (4.5, -0.5);
      \draw (3.5, -0.5) node[anchor=south] {$r_x$};
      
      \draw [<->] (-1.5, 0.5) -- (-1.5, 2.5);
      \draw (-1.5, 1.5) node[anchor=east] {$r_y$};
      
      \draw [dotted,color=gray] (2.5, -0.5) -- (2.5, 0.5);
      \draw [dotted,color=gray] (4.5, -0.5) -- (4.5, 0.5);
      
      \draw [dotted,color=gray] (0, -0.5) -- (0, 0);
      \draw [dotted,color=gray] (1, -0.5) -- (1, 0);
      
      \draw [dotted,color=gray] (-0.5, 0) -- (0, 0);
      \draw [dotted,color=gray] (-0.5, 1) -- (0, 1);
      
      \draw [dotted,color=gray] (-1.5, 0.5) -- (0.5, 0.5);
      \draw [dotted,color=gray] (-1.5, 2.5) -- (0.5, 2.5);
      
      \draw [dotted,color=gray] (9,0) -- (9.5, 0);
      \draw [dotted,color=gray] (9,7) -- (9.5, 7);
      \draw [dotted,color=gray] (9,7) -- (9, 7.5);
      \draw [dotted,color=gray] (0,7) -- (0, 7.5);
    \end{tikzpicture}
  \end{center}
  \caption{Exemplary sketch of the probe and input wave lattices drawn on top of the simulation domain.
           The probe positions $\mathcal{P}$ are denoted by circles and the input wave positions
           $\mathcal{I}$ by crosses. In this case, we have $P_x = 2R_x$ and $P_y = 2R_y$ as well as
           $\mathcal{I}_\text{offset} = (p_x/2, p_y/2)$. Note the
           compatibility of both lattices and that they cover the entire specimen domain, respecting the
           periodic boundaries.}\label{lma_lattices}
\end{figure}
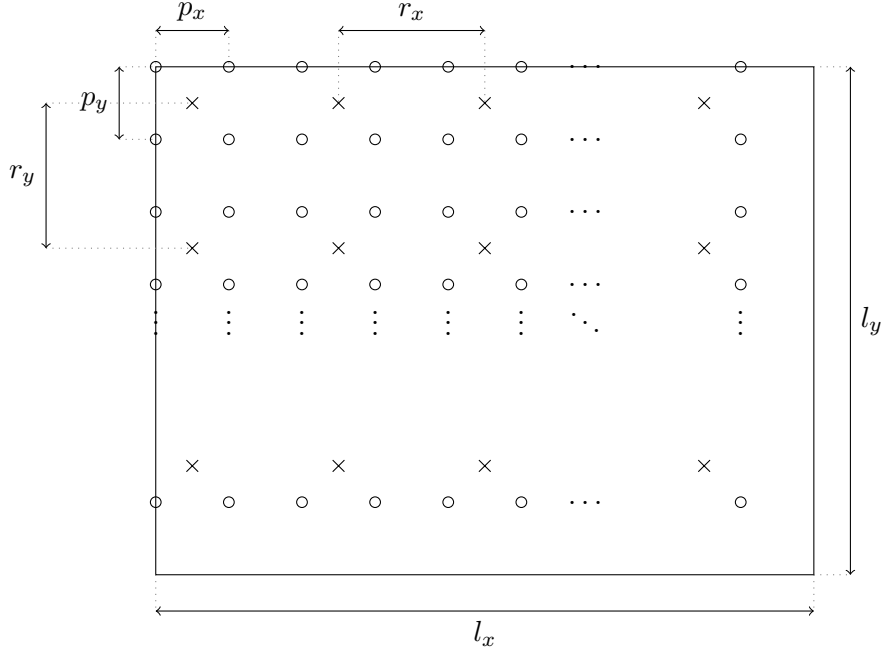

Although, from a mathematical point of view, the number of elements of $\mathcal{I}$ can
clearly be chosen smaller than or equal to the number of probe positions $\abs{\mathcal{P}}$
for a good approximation of the STEM probes, it is not clear that this is always beneficial
for the computations (see \cref{rsapprox_inputwaves} for a discussion of some possible input
waves). Therefore, we don't exclude the case $\abs{\mathcal{I}} > \abs{\mathcal{P}}$. The
probe and input wave lattices are schematically depicted in \cref{lma_lattices}.\newline

If we used a full real space basis
$(u_i)_{i\in\mathcal{I}}$ as in \cref{task_realspace_basis} for the approximation of the
STEM probes, then at most a minor improvement of the computation time compared to the standard Multislice algorithm is expected.
In PRISM, the computational speed up is achieved by reducing the number of basis elements from
$\abs{\mathcal{K}}$ to $\abs{\mathcal{K}_f}$, which accordingly reduces the number of times
that the Schrödinger equation needs to be solved. It is possible to trade off computation time
against accuracy in a similar way when working with localized real space functions $u_i$.
Among all the possibilities to do so, we note in particular the following three approaches:
\begin{enumerate}[(i)]
\item\label{lrs_opt1} Considering only every $f$-th point on the input wave lattice in both directions, thereby reducing the
                      number of evaluations of the Multislice algorithm by a factor of $f^2$.
\item\label{lrs_opt2} Performing the Multislice algorithm on a smaller computation window of the size $X'\times Y'$ for
                      $X',Y'\in\N$ with $X'<X$, $Y'<Y$, adjusted to the extent of the localized real space functions $u_i$.
\item\label{lrs_opt3} For a STEM probe $\hat\psi_p^\text{init}$, considering only the $L\in\N$ input waves $u_i$
                      closest to the probe position $p$ for the approximation of $\hat\psi_p^\text{init}$.
\end{enumerate}

The first option \ref{lrs_opt1} is very similar to PRISM's interpolation factor. The set
of the input wave positions attained by only considering every $f$-th point of the input
wave lattice is denoted by
\begin{equation*}
  \mathcal{I}_f := \left\{af\begin{pmatrix*} r_x \\ 0 \end{pmatrix*} + bf\begin{pmatrix*} 0 \\ r_y \end{pmatrix*} + \mathcal{I}_\text{offset} \;\middle|\; a\in\{0,\ldots,R_x/f-1\}, b\in\{0,\ldots,R_y/f-1\} \right\},
\end{equation*}
where we generally assume that $R_x$ and $R_y$ are multiples of $f$ by rounding up to the next multiple of $f$ if necessary and
adjusting the probe lattice accordingly.
Option \ref{lrs_opt2} is analogous to the approach of reducing the size of the
computation window by a factor of $f^2$; the only difference here is, that solutions to the
Schrödinger equation are computed for the input waves as initial condition instead
of the STEM probes. The third option \ref{lrs_opt3} is unique to this approach and only possible because both
the STEM probes and the input waves are localized in real space. It is motivated by the
idea that input waves close to the center of a STEM probe may contribute more to its
approximation than those further away. The
parameter $L$ can be freely chosen to adjust the number of input waves that
are considered for the approximation of a given STEM probe.

It should be noted that
\ref{lrs_opt1} and \ref{lrs_opt3} are not independent from each other. A larger
value of $f$ generally requires a larger value of $L$ to achieve the same approximation
quality and vice versa. The second option \ref{lrs_opt2} is mostly interesting in the case
where the Multislice algorithm is entirely performed in real space without use of the
Fourier transform, as the periodicity of the DFT likely does not match the periodicity of
the specimen for computation windows smaller than whole simulation domain. Nonetheless,
if $U_i$ denotes the solution to the Schrödinger equation with initial condition $u_i$,
then $U_i$ will also be fairly well localized in real space. This means that the result
can always be cropped to a smaller computation domain, which can save large amounts of
computer memory. While \ref{lrs_opt1} is the main
factor for reducing the computation times, \ref{lrs_opt3} also plays an important role
as can be seen in the expressions for the computational complexity below.\newline

All that has been described above leads to a new algorithm for the simulation of STEM images,
which is summarized in \cref{lma_pseudocode}. This algorithm is referred to as the Lattice
Multislice Algorithm (LMA) due to the important role of the probe and input wave lattices.

\begin{figure}[h]
  \hfill
  \begin{mdframed}[backgroundcolor=BoxBg, skipabove=10pt, skipbelow=10pt, leftmargin=5pt, rightmargin=5pt, roundcorner=5pt]
  \begin{minipage}{0.98\textwidth}
    Initialization: choose a localized real space function $u$, the input wave lattice $\mathcal{I}$ as
    well as $f\in\N$ and $L\in\N$. Precompute the approximation coefficients.
    \begin{enumerate}
      \item For each $i\in\mathcal{I}_f$, use the Multislice algorithm (cf. \cref{multislice_pseudocode})
            to calculate the solution $U_{i}$ to the Schrödinger equation with initial condition $u_{i}$.
      \item For each probe position $p\in \mathcal{P}$, calculate an approximation to the corresponding
            exit wave $\hat\psi_p^\text{exit}$ as a linear combination of the $U_i$ using the precomputed approximation coefficients.
    \end{enumerate}
  \end{minipage}
  \end{mdframed}
  \caption{Pseudocode for the Lattice Multislice Algorithm.}\label{lma_pseudocode}
\end{figure}

The exact formula for calculating the linear combination in the second step in \cref{lma_pseudocode} is given in
\cref{rsapprox_coefficients}.

\paragraph{Computational complexity} We keep the same notation as before, where the simulation domain is discretized to
$X \times Y$ pixel, and denote by $X' \times Y'$ the size of the smaller computation
domain for the input waves in pixel.
Disregarding the precomputation of the approximation coefficients, the number of necessary
floating point operations to compute one STEM image with the Lattice Multislice Algorithm
can be estimated as
\begin{align*}
  T_\text{LMA} &= \abs{\mathcal{I}_f}T_{\text{Multislice}, XY} + \abs{\mathcal{P}}LX'Y' \\
               &= \frac{1}{f^2}\abs{\mathcal{I}} N (2XY + 2XY\log(XY)) + \abs{\mathcal{P}}LX'Y'
\end{align*}
if the propagation step in the Multislice algorithm is performed in Fourier space using the FFT, or as
\begin{align*}
  \tilde T_\text{LMA} &= \abs{\mathcal{I}_f}\tilde T_{\text{Multislice}, X'Y'} + \abs{\mathcal{P}}LX'Y' \\
                      &= \frac{1}{f^2}\abs{\mathcal{I}} N (X'Y' + X'Y'K_1K_2) + \abs{\mathcal{P}}LX'Y'.
\end{align*}
if the propagation step is performed as a convolution in real space.

If we found a basis $\mathcal{I}$ with $\abs{\mathcal{K}}$ elements, chose $L=\abs{\mathcal{K}}/f^2$ and reduced
the computation domain to the same size as in PRISM, i.e. $X' = X/f$ and $Y' = Y/f$, then
\begin{align*}
  T_\text{LMA} &= \frac{1}{f^2}\abs{\mathcal{K}} N (2XY + 2XY\log(XY)) + \abs{\mathcal{P}}\frac{\abs{\mathcal{K}}}{f^2}\frac{X}{f}\frac{Y}{f} \\
               &= \frac{1}{f^2}\abs{\mathcal{K}} N (2XY + 2XY\log(XY)) + \frac{1}{f^4}\abs{\mathcal{P}}\abs{\mathcal{K}}XY \\
               &\approx T_\text{PRISM}.
\end{align*}

It is interesting to note that if we choose $X' = X/f$ and $Y' = Y/f$, then the first summand
of $\tilde T_\text{LMA}$ is proportional to $f^{-4}$. Comparing this to $T_\text{LMA}$
yields
\begin{align*}
  \tilde T_\text{LMA} &\le T_\text{LMA} \\
  \Longleftrightarrow\quad \frac{1}{f^2}\abs{\mathcal{I}} N (X'Y' + X'Y'K_1K_2) &\le \frac{1}{f^2}\abs{\mathcal{I}} N (2XY + 2XY\log(XY)) \\
  \Longleftrightarrow\quad \frac{1}{f^4}\abs{\mathcal{I}} N (XY + XYK_1K_2) &\le \frac{1}{f^2}\abs{\mathcal{I}} N (2XY + 2XY\log(XY)) \\
  \Longleftrightarrow\quad \frac{XY + XYK_1K_2}{2XY + 2XY\log(XY)} &\le f^2 \\
  \Longleftrightarrow\quad \sqrt{\frac{1}{2}\frac{1 + K_1K_2}{1 + 1\log(XY)}} &\le f,
\end{align*}
or, approximately,
\begin{equation*}
  f \gtrapprox \sqrt{\frac{K_1K_2}{2\log(XY)}}.
\end{equation*}
For a simulation domain size of $2048^2$ pixels and a discretized convolution kernel
with a width and height of $25$ pixels, the above inequality is true if and only if $f\ge 4$.

\subsubsection{Approximation coefficients}\label{rsapprox_coefficients}
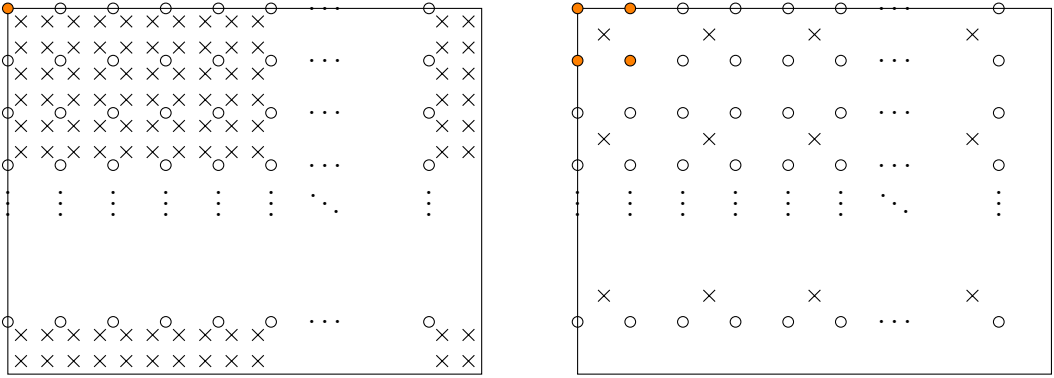
\begin{figure}[h]
  \begin{minipage}{0.49\textwidth}
  \begin{center}
    \begin{tikzpicture}[yscale=-1, x=1.8em, y=1.8em]
      \draw [draw=black] (0, 0) rectangle (9, 7);
      \draw [fill=orange] (0, 0) circle (2pt); \draw (1, 0) circle (2pt); \draw (2, 0) circle (2pt); \draw (3, 0) circle (2pt); \draw (4, 0) circle (2pt); \draw (5, 0) circle (2pt); \draw (8, 0) circle (2pt);
      \draw (0, 1) circle (2pt); \draw (1, 1) circle (2pt); \draw (2, 1) circle (2pt); \draw (3, 1) circle (2pt); \draw (4, 1) circle (2pt); \draw (5, 1) circle (2pt); \draw (8, 1) circle (2pt);
      \draw (0, 2) circle (2pt); \draw (1, 2) circle (2pt); \draw (2, 2) circle (2pt); \draw (3, 2) circle (2pt); \draw (4, 2) circle (2pt); \draw (5, 2) circle (2pt); \draw (8, 2) circle (2pt);
      \draw (0, 3) circle (2pt); \draw (1, 3) circle (2pt); \draw (2, 3) circle (2pt); \draw (3, 3) circle (2pt); \draw (4, 3) circle (2pt); \draw (5, 3) circle (2pt); \draw (8, 3) circle (2pt);
      \draw (0, 6) circle (2pt); \draw (1, 6) circle (2pt); \draw (2, 6) circle (2pt); \draw (3, 6) circle (2pt); \draw (4, 6) circle (2pt); \draw (5, 6) circle (2pt); \draw (8, 6) circle (2pt);
      
      \draw (5.5, 0) node[anchor=west] {$\ldots$};
      \draw (5.5, 1) node[anchor=west] {$\ldots$};
      \draw (5.5, 2) node[anchor=west] {$\ldots$};
      \draw (5.5, 3) node[anchor=west] {$\ldots$};
      \draw (5.5, 6) node[anchor=west] {$\ldots$};
      
      \draw (0, 3) node[anchor=north] {$\vdots$};
      \draw (1, 3) node[anchor=north] {$\vdots$};
      \draw (2, 3) node[anchor=north] {$\vdots$};
      \draw (3, 3) node[anchor=north] {$\vdots$};
      \draw (4, 3) node[anchor=north] {$\vdots$};
      \draw (5, 3) node[anchor=north] {$\vdots$};
      \draw (8, 3) node[anchor=north] {$\vdots$};
      
      \draw (5.5, 3) node[anchor=north west] {$\ddots$};
      
      \draw (0.25, 0.25) node[cross] {}; \draw (0.75, 0.25) node[cross] {}; \draw (1.25, 0.25) node[cross] {}; \draw (1.75, 0.25) node[cross] {}; \draw (2.25, 0.25) node[cross] {}; \draw (2.75, 0.25) node[cross] {}; \draw (3.25, 0.25) node[cross] {}; \draw (3.75, 0.25) node[cross] {}; \draw (4.25, 0.25) node[cross] {}; \draw (4.75, 0.25) node[cross] {}; \draw (8.25, 0.25) node[cross] {}; \draw (8.75, 0.25) node[cross] {};
      \draw (0.25, 0.75) node[cross] {}; \draw (0.75, 0.75) node[cross] {}; \draw (1.25, 0.75) node[cross] {}; \draw (1.75, 0.75) node[cross] {}; \draw (2.25, 0.75) node[cross] {}; \draw (2.75, 0.75) node[cross] {}; \draw (3.25, 0.75) node[cross] {}; \draw (3.75, 0.75) node[cross] {}; \draw (4.25, 0.75) node[cross] {}; \draw (4.75, 0.75) node[cross] {}; \draw (8.25, 0.75) node[cross] {}; \draw (8.75, 0.75) node[cross] {};
      \draw (0.25, 1.25) node[cross] {}; \draw (0.75, 1.25) node[cross] {}; \draw (1.25, 1.25) node[cross] {}; \draw (1.75, 1.25) node[cross] {}; \draw (2.25, 1.25) node[cross] {}; \draw (2.75, 1.25) node[cross] {}; \draw (3.25, 1.25) node[cross] {}; \draw (3.75, 1.25) node[cross] {}; \draw (4.25, 1.25) node[cross] {}; \draw (4.75, 1.25) node[cross] {}; \draw (8.25, 1.25) node[cross] {}; \draw (8.75, 1.25) node[cross] {};
      \draw (0.25, 1.75) node[cross] {}; \draw (0.75, 1.75) node[cross] {}; \draw (1.25, 1.75) node[cross] {}; \draw (1.75, 1.75) node[cross] {}; \draw (2.25, 1.75) node[cross] {}; \draw (2.75, 1.75) node[cross] {}; \draw (3.25, 1.75) node[cross] {}; \draw (3.75, 1.75) node[cross] {}; \draw (4.25, 1.75) node[cross] {}; \draw (4.75, 1.75) node[cross] {}; \draw (8.25, 1.75) node[cross] {}; \draw (8.75, 1.75) node[cross] {};
      \draw (0.25, 2.25) node[cross] {}; \draw (0.75, 2.25) node[cross] {}; \draw (1.25, 2.25) node[cross] {}; \draw (1.75, 2.25) node[cross] {}; \draw (2.25, 2.25) node[cross] {}; \draw (2.75, 2.25) node[cross] {}; \draw (3.25, 2.25) node[cross] {}; \draw (3.75, 2.25) node[cross] {}; \draw (4.25, 2.25) node[cross] {}; \draw (4.75, 2.25) node[cross] {}; \draw (8.25, 2.25) node[cross] {}; \draw (8.75, 2.25) node[cross] {};
      \draw (0.25, 2.75) node[cross] {}; \draw (0.75, 2.75) node[cross] {}; \draw (1.25, 2.75) node[cross] {}; \draw (1.75, 2.75) node[cross] {}; \draw (2.25, 2.75) node[cross] {}; \draw (2.75, 2.75) node[cross] {}; \draw (3.25, 2.75) node[cross] {}; \draw (3.75, 2.75) node[cross] {}; \draw (4.25, 2.75) node[cross] {}; \draw (4.75, 2.75) node[cross] {}; \draw (8.25, 2.75) node[cross] {}; \draw (8.75, 2.75) node[cross] {};
      
      \draw (0.25, 6.25) node[cross] {}; \draw (0.75, 6.25) node[cross] {}; \draw (1.25, 6.25) node[cross] {}; \draw (1.75, 6.25) node[cross] {}; \draw (2.25, 6.25) node[cross] {}; \draw (2.75, 6.25) node[cross] {}; \draw (3.25, 6.25) node[cross] {}; \draw (3.75, 6.25) node[cross] {}; \draw (4.25, 6.25) node[cross] {}; \draw (4.75, 6.25) node[cross] {}; \draw (8.25, 6.25) node[cross] {}; \draw (8.75, 6.25) node[cross] {};
      \draw (0.25, 6.75) node[cross] {}; \draw (0.75, 6.75) node[cross] {}; \draw (1.25, 6.75) node[cross] {}; \draw (1.75, 6.75) node[cross] {}; \draw (2.25, 6.75) node[cross] {}; \draw (2.75, 6.75) node[cross] {}; \draw (3.25, 6.75) node[cross] {}; \draw (3.75, 6.75) node[cross] {}; \draw (4.25, 6.75) node[cross] {}; \draw (4.75, 6.75) node[cross] {}; \draw (8.25, 6.75) node[cross] {}; \draw (8.75, 6.75) node[cross] {};
      
    \end{tikzpicture}
  \end{center}
  \end{minipage}
  \hfill
  \begin{minipage}{0.49\textwidth}
  \begin{center}
    \begin{tikzpicture}[yscale=-1, x=1.8em, y=1.8em]
      \draw [draw=black] (0, 0) rectangle (9, 7);
      \draw [fill=orange] (0, 0) circle (2pt); \draw [fill=orange] (1, 0) circle (2pt); \draw (2, 0) circle (2pt); \draw (3, 0) circle (2pt); \draw (4, 0) circle (2pt); \draw (5, 0) circle (2pt); \draw (8, 0) circle (2pt);
      \draw [fill=orange] (0, 1) circle (2pt); \draw [fill=orange] (1, 1) circle (2pt); \draw (2, 1) circle (2pt); \draw (3, 1) circle (2pt); \draw (4, 1) circle (2pt); \draw (5, 1) circle (2pt); \draw (8, 1) circle (2pt);
      \draw (0, 2) circle (2pt); \draw (1, 2) circle (2pt); \draw (2, 2) circle (2pt); \draw (3, 2) circle (2pt); \draw (4, 2) circle (2pt); \draw (5, 2) circle (2pt); \draw (8, 2) circle (2pt);
      \draw (0, 3) circle (2pt); \draw (1, 3) circle (2pt); \draw (2, 3) circle (2pt); \draw (3, 3) circle (2pt); \draw (4, 3) circle (2pt); \draw (5, 3) circle (2pt); \draw (8, 3) circle (2pt);
      \draw (0, 6) circle (2pt); \draw (1, 6) circle (2pt); \draw (2, 6) circle (2pt); \draw (3, 6) circle (2pt); \draw (4, 6) circle (2pt); \draw (5, 6) circle (2pt); \draw (8, 6) circle (2pt);
      
      \draw (5.5, 0) node[anchor=west] {$\ldots$};
      \draw (5.5, 1) node[anchor=west] {$\ldots$};
      \draw (5.5, 2) node[anchor=west] {$\ldots$};
      \draw (5.5, 3) node[anchor=west] {$\ldots$};
      \draw (5.5, 6) node[anchor=west] {$\ldots$};
      
      \draw (0, 3) node[anchor=north] {$\vdots$};
      \draw (1, 3) node[anchor=north] {$\vdots$};
      \draw (2, 3) node[anchor=north] {$\vdots$};
      \draw (3, 3) node[anchor=north] {$\vdots$};
      \draw (4, 3) node[anchor=north] {$\vdots$};
      \draw (5, 3) node[anchor=north] {$\vdots$};
      \draw (8, 3) node[anchor=north] {$\vdots$};
      
      \draw (5.5, 3) node[anchor=north west] {$\ddots$};
      
      \draw (0.5, 0.5) node[cross] {}; \draw (2.5, 0.5) node[cross] {}; \draw (4.5, 0.5) node[cross] {}; \draw (7.5, 0.5) node[cross] {};
      \draw (0.5, 2.5) node[cross] {}; \draw (2.5, 2.5) node[cross] {}; \draw (4.5, 2.5) node[cross] {}; \draw (7.5, 2.5) node[cross] {};
      \draw (0.5, 5.5) node[cross] {}; \draw (2.5, 5.5) node[cross] {}; \draw (4.5, 5.5) node[cross] {}; \draw (7.5, 5.5) node[cross] {};
    \end{tikzpicture}
  \end{center}
  \end{minipage}
  \caption{Example for \bref{lma_approximation_case1} on the left and \bref{lma_approximation_case2} on the right. Circles
           indicate probe positions and crosses indicate input wave positions in $\mathcal{I}_f$ drawn on top of the
           simulation domain depicted by the large rectangle.
           The probe positions with unique sets of coefficients as in
           \cref{lma_approximation_exit1_0,lma_approximation_exit2_0} are denoted by the filled circles.
           For the sake of simplicity, we disregard rotational and mirror symmetries in the coefficients,
           which would reduce the required amount of precomputations even further.}\label{lma_lattices_approximation}
\end{figure}

In principle, it is possible to calculate approximation coefficients $\alpha_{p, i}\in\C$
for all probe positions $p\in\mathcal{P}$ and all $i\in\mathcal{I}_f$ such that
\begin{equation*}
  \hat\psi_p^\text{init} \approx \sum_{i\in\mathcal{I}_{f}} \alpha_{p, i} u_i
\end{equation*}
for all $p\in\mathcal{P}$ using linear least squares. However, this would be a waste of computations
due to the compatibility property of the probe and input wave lattices, which makes it possible to
reuse the same approximation coefficients for many probe positions. To show how the coefficients can
be reused, we distinguish
\bref{lma_approximation_case1} with $\abs{\mathcal{I}} \ge \abs{\mathcal{P}}$ and \bref{lma_approximation_case2}
with $\abs{\mathcal{P}} \ge \abs{\mathcal{I}}$. Both cases are illustrated in \cref{lma_lattices_approximation}
and described below.

Denote by $\mathcal{I}_{f, p, L}$ for $p\in\mathcal{P}$ the subset
of $\mathcal{I}_f$ containing the $L$ closest input waves to $p$.
Looking at \bref{lma_approximation_case1} in more detail, it is easy to see that only the coefficients
for the approximation of a single probe need to be calculated and can in turn be reused for all other
probe positions. Explicitly, if
\begin{equation}\label{lma_approximation_exit1_0}
  \hat\psi_0^\text{init} \approx \sum_{i\in\mathcal{I}_{f, 0, L}} \alpha_i u_i
\end{equation}
for suitable $\alpha_i\in\C$, then
\begin{equation}\label{lma_approximation_exit1}
  \hat\psi_p^\text{init} \approx \sum_{i\in\mathcal{I}_{f, 0, L}} \alpha_i u_{i+p}.
\end{equation}
with precisely the same approximation error for all $p\in\mathcal{P}$. Here, $i+p$ is
understood to be calculated in $\R^2 / (l_x\Z \times l_y\Z)$ and again an element of
$\mathcal{I}_f$ due to the compatibility of the lattices and \bref{lma_lattices_cover_specimen}.

For \bref{lma_approximation_case2}, in total $c \cdot d$ unique sets of coefficients need to be calculated for the probe positions
in
\begin{equation*}
  \mathcal{P}_\text{approx} := \left\{a\begin{pmatrix*} p_x \\ 0 \end{pmatrix*} + b\begin{pmatrix*} 0 \\ p_y \end{pmatrix*} \;\middle|\; a\in\{0,\ldots,c-1\}, \; b\in\{0,\ldots,d-1\} \right\}.
\end{equation*}
If the approximations are given by
\begin{equation}\label{lma_approximation_exit2_0}
  \hat\psi_p^\text{init} \approx \sum_{i\in\mathcal{I}_{f, p, L}} \alpha_{p, i} u_i \qquad\forall\, p\in\mathcal{P}_\text{approx}
\end{equation}
for suitable $\alpha_{p, i}\in\C$, then, similarly to \cref{lma_approximation_exit1},
\begin{equation}\label{lma_approximation_exit2}
  \hat\psi_{p+q}^\text{init} \approx \sum_{i\in\mathcal{I}_{f, p, L}} \alpha_{p, i} u_{i+q} \qquad\forall\, p\in\mathcal{P}_\text{approx}
\end{equation}
holds for all $q = xcp_xe_1+ydp_ye_2$ with $x, y\in\N$ such that $p+q\in\mathcal{P}$, where $e_1 = (1\; 0)$ and $e_2 = (0\; 1)$.
As above, $i+q$ is understood to be calculated as an element of $\R^2 / (l_x\Z \times l_y\Z)$ and
again contained in $\mathcal{I}_f$.

\subsubsection{Input wave functions}\label{rsapprox_inputwaves}
The description of the Lattice Multislice Algorithm is incomplete without a discussion of
possible choices for the input waves. By definition, the STEM probes $\hat\psi_p^\text{init}$
form a generating system for $V_\mathcal{P}$, so an obvious choice is $u = \hat\psi_0^\text{init}$.
This is an interesting choice because it makes the approximation very easy by requiring only
comparatively small values of $L$ as is shown below.
If we additionally choose $\mathcal{I} = \mathcal{P}$,
$f=1$ and $L=1$, then we retrieve the standard algorithm for the simulation of STEM images,
where the Schrödinger equation is solved once for every probe position with the STEM
probes itself as initial conditions.

A more interesting case arises if we deviate from these
parameters, in particular if we increase $f$ (and, necessarily, $L$). In order to avoid
a large spread in the approximation errors for different probe positions, we choose the input
wave lattice to be identical to the probe lattice, but shifted by half of the probe spacing $p_x/2$ and
$p_y/2$ in each direction respectively. Explicitly, we set $\mathcal{I} := \mathcal{P} + \frac{1}{2}(p_x, p_y)$.
Without this shift, the approximation would be perfect for every probe position $p\in\mathcal{P}\cap\mathcal{I}_f$,
but subject to an error for all probe positions $p\in\mathcal{P}\backslash\mathcal{I}_f$. This may cause
potentially significant artifacts in a checkerboard pattern in the simulated STEM image for $f>1$.

\begin{figure}
  \begin{center}
    \begin{tikzpicture}
      \begin{axis}[width = 0.95\textwidth, height = 16em,
                   legend pos = outer north east,
                   ymode = log,
                   xmin = 0, xmax = 100,
                   ymin = 8e-4, ymax = 1]
        \addplot[mark=none, color=orange] table[x expr=\thisrowno{0}, y expr=\thisrowno{1}] {dat_AT_PSf1_approximation_errors_euc};
        \addlegendentry{$f=1$}
        \addplot[mark=none, color=blue] table[x expr=\thisrowno{0}, y expr=\thisrowno{1}] {dat_AT_PSf2_approximation_errors_euc};
        \addlegendentry{$f=2$}
        \addplot[mark=none, color=green] table[x expr=\thisrowno{0}, y expr=\thisrowno{1}] {dat_AT_PSf3_approximation_errors_euc};
        \addlegendentry{$f=3$}
        \addplot[mark=none, color=black] table[x expr=\thisrowno{0}, y expr=\thisrowno{1}] {dat_AT_PSf4_approximation_errors_euc};
        \addlegendentry{$f=4$}
      \end{axis}
    \end{tikzpicture}
    \begin{tikzpicture}
      \begin{axis}[width = 0.95\textwidth, height = 16em,
                   legend pos = outer north east,
                   ymode = log,
                   xmin = 0, xmax = 100,
                   ymin = 1e-3, ymax = 1]
        \addplot[mark=none, color=orange] table[x expr=\thisrowno{0}, y expr=\thisrowno{1}] {dat_AT_PSf1_approximation_errors_sup};
        \addlegendentry{$f=1$}
        \addplot[mark=none, color=blue] table[x expr=\thisrowno{0}, y expr=\thisrowno{1}] {dat_AT_PSf2_approximation_errors_sup};
        \addlegendentry{$f=2$}
        \addplot[mark=none, color=green] table[x expr=\thisrowno{0}, y expr=\thisrowno{1}] {dat_AT_PSf3_approximation_errors_sup};
        \addlegendentry{$f=3$}
        \addplot[mark=none, color=black] table[x expr=\thisrowno{0}, y expr=\thisrowno{1}] {dat_AT_PSf4_approximation_errors_sup};
        \addlegendentry{$f=4$}
      \end{axis}
    \end{tikzpicture}
  \end{center}
  \caption{Relative approximation errors for the approximation of the STEM probe $\hat\psi_0^\text{init}$
           in the euclidean norm (top) and the supremum norm (bottom) for various values of
           $f$ and $L$. The x axis depicts the value of $L$ and the y axis the approximation
           error. The approximations have been computed with linear least squares. The parameters used are $u = \hat\psi_0^\text{init}$ and
           $\mathcal{I} = \mathcal{P} + \frac{1}{2}(p_x, p_y)$ with $X=Y=1024$ and
           $P_x = P_y = 524$. The simulation window size is $l_x = l_y = 62.416$ (in Angstrom)
           and the microscope settings are $\lambda = 0.0250793$ (corresponding to an accelerating
           voltage of $U = 200000$ Volt), $Z = 100$, $C_s = -2000$ and $\alpha_\text{max} = 0.026$.}\label{lma_inputwave_probe_approximation}
\end{figure}
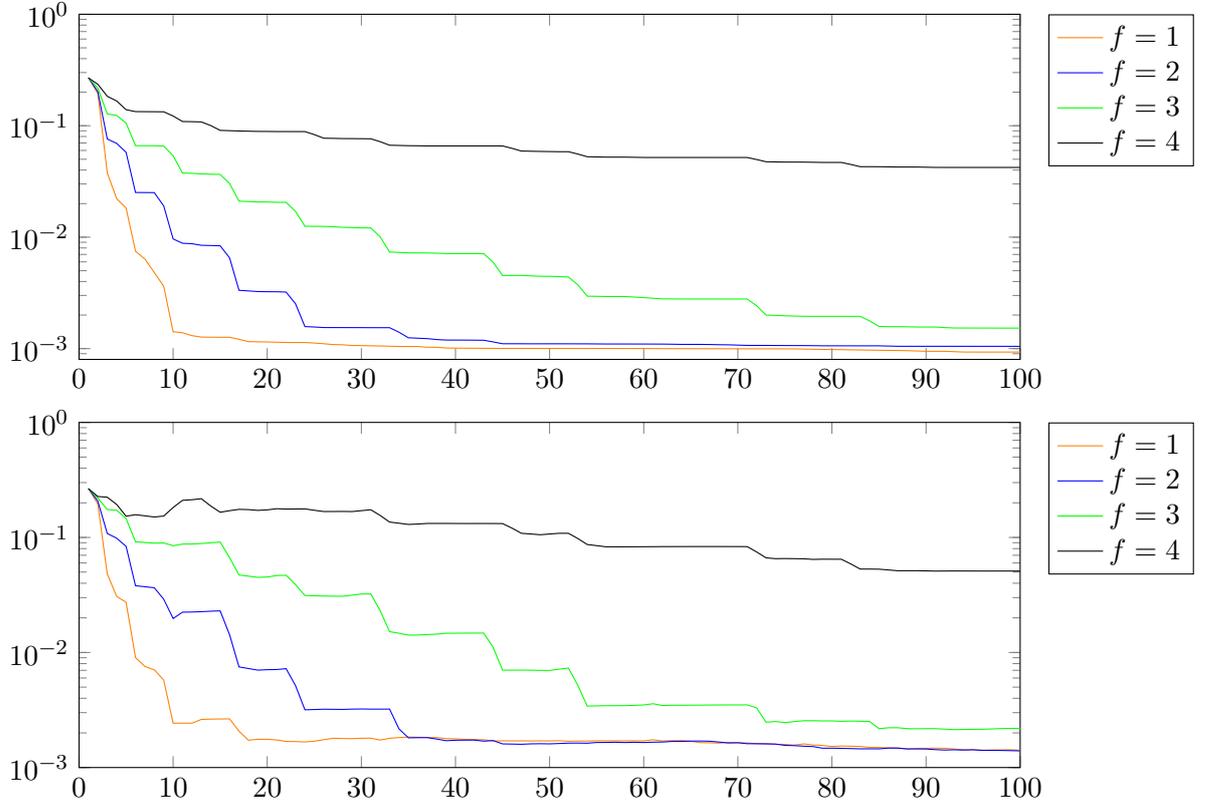

An example of the magnitude of the approximation error for one particular kind of parameters
for $u = \hat\psi_0^\text{init}$
is shown in \cref{lma_inputwave_probe_approximation} for varying values of $f$ and $L$. When increasing
$L$, the approximation error generally appears to decrease in a staircase fashion: increasing
$L$ by 1 either results in a minimal or a very large reduction of the approximation error.

Input waves with smaller support or faster decay to zero in real space than the STEM
probes require less computer memory and accelerate the computations if the propagation
step in the Multislice algorithm is performed as a convolution in real space. However, 
input waves with a smaller support will in general require larger values of $L$ or smaller
values of $f$ in order to achieve the same approximation error as other input waves with
larger values of $f$ or smaller $L$. This in turn increases the required computer memory
and computation time. It is unclear what the optimal choice of $u$ is, both with respect
to the use of computer memory and computation time.

We present a few nontrivial choices
of input waves in the following, alongside with numerical experiments depicting the
tradeoff of smaller support vs. the choice of of $L$ and $f$ and the resulting approximation
errors.

\paragraph{Trigonometric polynomials} The trigonometric polynomials
introduced in \cite{rauhut05} are an attractive choice for the input waves due to their
property of being optimally localized in real space for a given degree. They may be
extended to a two-dimensional function in several ways. One option is to define $u$
as the tensor product $u(x, y) := \varphi_n(x)\varphi_n(y)$, where
\begin{equation*}
  \varphi_n(t) = \sum_{k=-n}^n e^{ikt}\cos\left(\frac{k\pi}{2n+2}\right) \qquad\forall\,t\in\R
\end{equation*}
is the trigonometric polynomial of degree $n\in\N$ defined in \cite{rauhut05}. This option has
the advantage that the approximation coefficients can be computed directly from the formulas for
the frequencies of the trigonometric polynomials and the STEM probes, not necessitating a least
squares approximation. Another option
is to define $u$ as the circularly symmetric function given by $u(v) := \varphi_n(\norm{v}_2)$.
The degree $n$ of $\varphi_n$ determines its maximum frequency and can be
calculated from the maximum frequency of the STEM probes and the Fourier space pixel size
$q\in\R$ as
\begin{equation*}
  n = \frac{\alpha_\text{max}}{\lambda} \cdot \frac{1}{q}.
\end{equation*}
This way, the frequency range of the input waves matches the
frequency range of the STEM probes optimally, ensuring a good approximation without
wasting computations on unneeded frequencies.

\begin{figure}
  \begin{center}
    \includegraphics[width=0.4\textwidth, keepaspectratio]{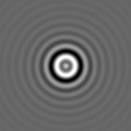}
    \quad
    \includegraphics[width=0.4\textwidth, keepaspectratio]{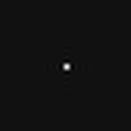}
  \end{center}
  \caption{The real part of the STEM probe in real space (left) and the input wave
           $u(x, y) := \varphi_n(x)\varphi_n(y)$ with $n=7$. Both images are
           discretized to $131 \times 131$ pixels with a pixel size of
           $0.1219062 \times 0.1219062$ (in Angstrom) and cover a similar frequency range.
           The parameters used for the computation of the STEM probe are
           $\lambda = 0.0250793$, $Z = 100$, $C_s = -2000$ and $\alpha_\text{max} = 0.026$.}
  \label{trigpoly_support_comparison}
\end{figure}

\cref{trigpoly_support_comparison} shows an image of a STEM probe next to the tensor
product of a trigono\-metric polynomial as described above with a matching frequency range.
Since this type of input wave is much better localized in real space than the probe,
it is important to carefully choose the probe and input wave lattices as well as
$f$ such the number of input waves close to a probe position is neither too small nor too large.
If the number is too small, then the input waves are spaced too far apart, resulting in a bad
approximation. If the number is too large, computation time is wasted. Another helpful perspective is
provided by the Fourier space representation: after discretization, the STEM probes are composed
of $N\in\N$ discrete frequencies within the aperture radius. The probe lattice and input wave lattice
as well as $f$ should be chosen such that the number of input waves is close to $N$ within the area around
the probe where the approximation is calculated.

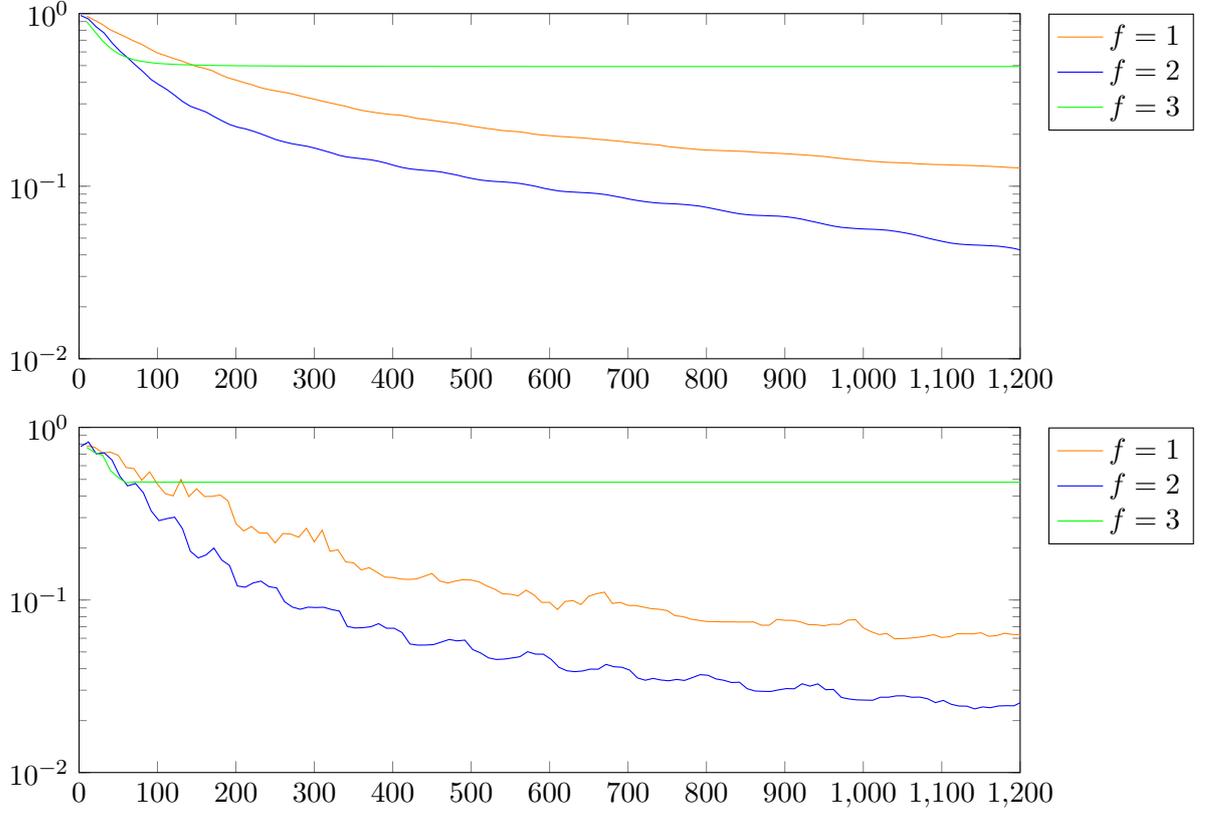
\begin{figure}
  \begin{center}
    \begin{tikzpicture}
      \begin{axis}[width = 0.95\textwidth, height = 16em,
                   legend pos = outer north east,
                   ymode = log,
                   xmin = 0, xmax = 1200,
                   ymin = 1e-2, ymax = 1]
        \addplot[mark=none, color=orange] table[x expr=\thisrowno{0}, y expr=\thisrowno{1}] {dat_AT_TPf1_approximation_errors_euc};
        \addlegendentry{$f=1$}
        \addplot[mark=none, color=blue] table[x expr=\thisrowno{0}, y expr=\thisrowno{1}] {dat_AT_TPf2_approximation_errors_euc};
        \addlegendentry{$f=2$}
        \addplot[mark=none, color=green] table[x expr=\thisrowno{0}, y expr=\thisrowno{1}] {dat_AT_TPf3_approximation_errors_euc};
        \addlegendentry{$f=3$}
      \end{axis}
    \end{tikzpicture}
    \begin{tikzpicture}
      \begin{axis}[width = 0.95\textwidth, height = 16em,
                   legend pos = outer north east,
                   ymode = log,
                   xmin = 0, xmax = 1200,
                   ymin = 1e-2, ymax = 1]
        \addplot[mark=none, color=orange] table[x expr=\thisrowno{0}, y expr=\thisrowno{1}] {dat_AT_TPf1_approximation_errors_sup};
        \addlegendentry{$f=1$}
        \addplot[mark=none, color=blue] table[x expr=\thisrowno{0}, y expr=\thisrowno{1}] {dat_AT_TPf2_approximation_errors_sup};
        \addlegendentry{$f=2$}
        \addplot[mark=none, color=green] table[x expr=\thisrowno{0}, y expr=\thisrowno{1}] {dat_AT_TPf3_approximation_errors_sup};
        \addlegendentry{$f=3$}
      \end{axis}
    \end{tikzpicture}
  \end{center}
  \caption{Relative approximation errors for the approximation of the STEM probe $\hat\psi_0^\text{init}$
           in the euclidean norm (top) and the supremum norm (bottom) for various values of
           $f$ and $L$. The x axis depicts the value of $L$ and the y axis the approximation
           error. The approximations have been computed with linear least squares. The parameters
           used are $u(x, y) = \varphi_n(x)\varphi_n(y)$ with $n=7$ and
           $\mathcal{I} = \mathcal{P} + \frac{1}{2}(p_x, p_y)$ with $X=Y=1024$ and
           $P_x = P_y = 314$. The simulation window size is $l_x = l_y = 62.416$ (in Angstrom)
           and the microscope settings are $\lambda = 0.0250793$ (corresponding to an accelerating
           voltage of $U = 200000$ Volt), $Z = 100$, $C_s = -2000$ and $\alpha_\text{max} = 0.026$.}\label{lma_inputwave_trigpoly_approximation}
\end{figure}

The consequences of this effect can be seen in
\cref{lma_inputwave_trigpoly_approximation}, where the optimal spacing in between the input
waves is achieved for $f=2$. This also shows that simply choosing the $L\in\N$ input waves closest to
the probe position as described in \cref{lrs_opt3} in \cref{lma_description} may yield bad approximation results. However, it
is unclear how a better subset of the input waves for a given configuration of probe lattice,
input wave lattice and input wave type can be found.

\paragraph{Gaussians}

When looking for localized functions in real space, it is natural to consider
Gaussians
\begin{equation*}
  g_\sigma(x) := \exp\left(\frac{-\norm{x}_2^2}{2\sigma^2}\right) \qquad\forall\,x\in\R^2
\end{equation*}
with variance $\sigma>0$.
Although technically non-zero everywhere, in practical computations they
can be considered equal to zero sufficiently far away from their center. Gaussians have
the useful property that their Fourier transform is again a Gaussian; similarly to the
trigonometric polynomials, this makes it possible to choose a reasonable value for the
variance based on the frequency range of the STEM probes, for example
\begin{equation*}
  \sigma = \frac{\lambda}{2\alpha_\text{max}}.
\end{equation*}

\begin{figure}
  \begin{center}
    \begin{tikzpicture}
      \begin{axis}[width = 0.95\textwidth, height = 16em,
                   legend pos = outer north east,
                   ymode = log,
                   xmin = 0, xmax = 1200,
                   ymin = 1e-2, ymax = 1]
        \addplot[mark=none, color=orange] table[x expr=\thisrowno{0}, y expr=\thisrowno{1}] {dat_AT_Gf1_approximation_errors_euc};
        \addlegendentry{$f=1$}
        \addplot[mark=none, color=blue] table[x expr=\thisrowno{0}, y expr=\thisrowno{1}] {dat_AT_Gf2_approximation_errors_euc};
        \addlegendentry{$f=2$}
        \addplot[mark=none, color=green] table[x expr=\thisrowno{0}, y expr=\thisrowno{1}] {dat_AT_Gf3_approximation_errors_euc};
        \addlegendentry{$f=3$}
        \addplot[mark=none, color=black] table[x expr=\thisrowno{0}, y expr=\thisrowno{1}] {dat_AT_Gf4_approximation_errors_euc};
        \addlegendentry{$f=4$}
      \end{axis}
    \end{tikzpicture}
    \begin{tikzpicture}
      \begin{axis}[width = 0.95\textwidth, height = 16em,
                   legend pos = outer north east,
                   ymode = log,
                   xmin = 0, xmax = 1200,
                   ymin = 1e-2, ymax = 1]
        \addplot[mark=none, color=orange] table[x expr=\thisrowno{0}, y expr=\thisrowno{1}] {dat_AT_Gf1_approximation_errors_sup};
        \addlegendentry{$f=1$}
        \addplot[mark=none, color=blue] table[x expr=\thisrowno{0}, y expr=\thisrowno{1}] {dat_AT_Gf2_approximation_errors_sup};
        \addlegendentry{$f=2$}
        \addplot[mark=none, color=green] table[x expr=\thisrowno{0}, y expr=\thisrowno{1}] {dat_AT_Gf3_approximation_errors_sup};
        \addlegendentry{$f=3$}
        \addplot[mark=none, color=black] table[x expr=\thisrowno{0}, y expr=\thisrowno{1}] {dat_AT_Gf4_approximation_errors_sup};
        \addlegendentry{$f=4$}
      \end{axis}
    \end{tikzpicture}
  \end{center}
  \caption{Relative approximation errors for the approximation of the STEM probe $\hat\psi_0^\text{init}$
           in the euclidean norm (top) and the supremum norm (bottom) for various values of
           $f$ and $L$. The x axis depicts the value of $L$ and the y axis the approximation
           error. The approximations have been computed with linear least squares. The parameters
           used are $u = g_\sigma$ with $\sigma=0.482$ Angstrom and
           $\mathcal{I} = \mathcal{P} + \frac{1}{2}(p_x, p_y)$ with $X=Y=1024$ and
           $P_x = P_y = 314$. The simulation window size is $l_x = l_y = 62.416$ (in Angstrom)
           and the microscope settings are $\lambda = 0.0250793$ (corresponding to an accelerating
           voltage of $U = 200000$ Volt), $Z = 100$, $C_s = -2000$ and $\alpha_\text{max} = 0.026$.}\label{lma_inputwave_gaussian_approximation}
\end{figure}
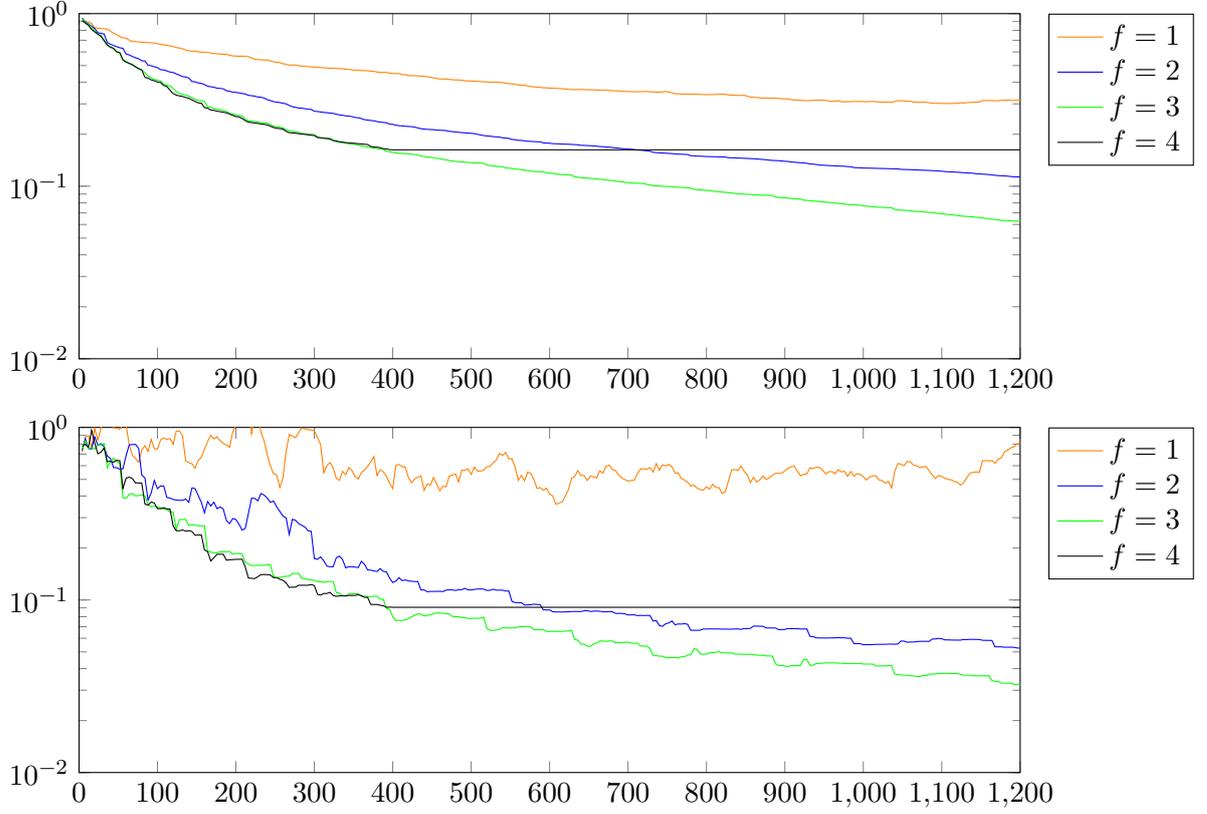

The approximation results for various values of $f$ and $L$ are plotted in
\cref{lma_inputwave_gaussian_approximation}.

\paragraph{Pixel input wave}

The fact that practically any function can be chosen as an input wave makes the Lattice Multislice
Algorithm highly versatile: the STEM probes themselves, trigonometric polynomials, Gaussians, indicator functions
of any set and other functions can be chosen as input waves. Since it is difficult to pinpoint the optimal choice of input waves for given STEM probe
parameters and given targets such as minimal memory consumption or computation time, it is interesting to look at the
extreme cases. Two of these were already discussed: (1) using the STEM probes themselves as input waves with
matching probe and input wave lattices $\mathcal{I} = \mathcal{P}$ with $f=1$ and $L=1$ and (2) using Fourier space
Dirac deltas as input waves as PRISM does. In the first case, we retrieve the original algorithm for STEM image
simulation where no approximation as in PRISM or LMA. The second case corresponds to input waves that are optimally
localized in Fourier space for a given discretization.

It is equally possible to consider Dirac deltas in real space as input waves, which then correspond to optimally
localized input waves in real space for a given discretization. This makes it possible to perform the Multislice
computations on a small subset of the simulation window, reducing the amount of computations for every individual
Multislice computation. However, this type of input wave requires one input wave for every pixel in the simulation
window. Although it is trivial to form the STEM probes from linear combinations of their pixel values and this
results in a perfect approximation, this type of input wave is impractical because of the sheer number of
necessary Multislice computations and summands in the linear combinations.


\paragraph{Probe approximation error and error in the output image}

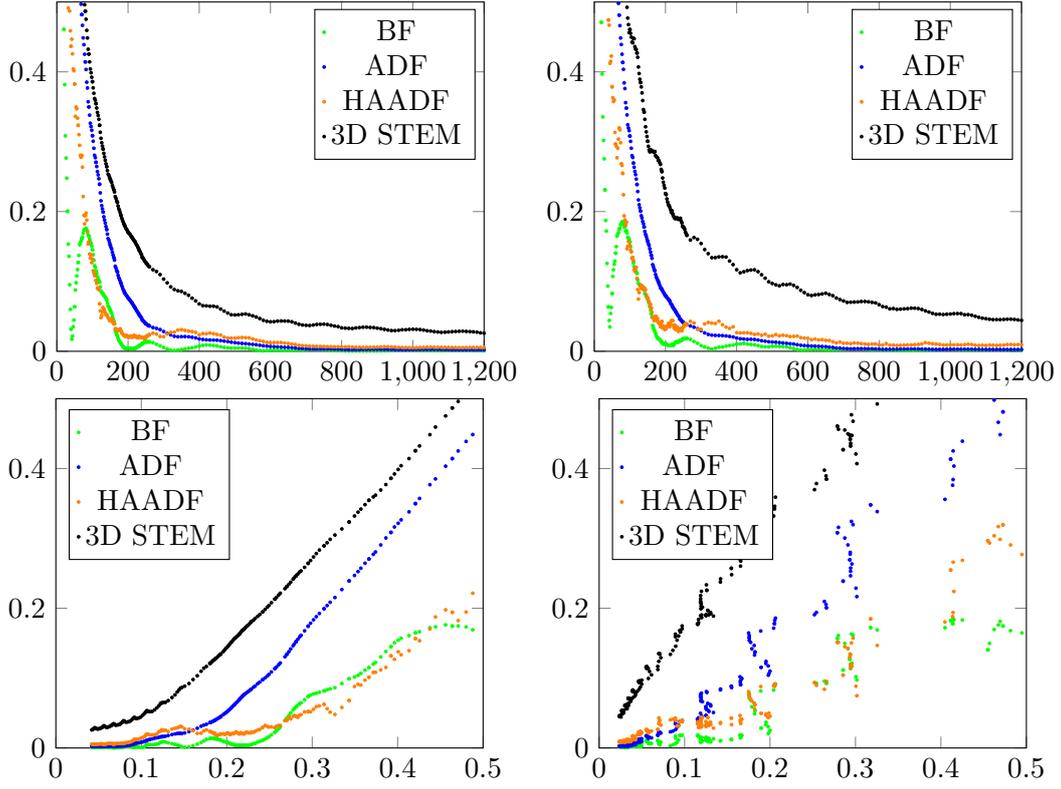
\begin{figure}
  \begin{center}
    \begin{tikzpicture}
      \begin{axis}[width = 0.49\textwidth,
                   xmin = 0, xmax = 1200,
                   ymin = 0, ymax = 0.5]
        \addplot[only marks, mark size = 0.5pt, color=green] table[x expr=\thisrowno{0}, y expr=\thisrowno{1}] {dat_AT_TPM1f2_bf_euc};
        \addlegendentry{BF}
        \addplot[only marks, mark size = 0.5pt, color=blue] table[x expr=\thisrowno{0}, y expr=\thisrowno{1}] {dat_AT_TPM1f2_adf_euc};
        \addlegendentry{ADF}
        \addplot[only marks, mark size = 0.5pt, color=orange] table[x expr=\thisrowno{0}, y expr=\thisrowno{1}] {dat_AT_TPM1f2_haadf_euc};
        \addlegendentry{HAADF}
        \addplot[only marks, mark size = 0.5pt, color=black] table[x expr=\thisrowno{0}, y expr=\thisrowno{1}] {dat_AT_TPM1f2_3d_euc};
        \addlegendentry{3D STEM}
        \addplot[only marks, mark size = 0.5pt, color=green] table[x expr=\thisrowno{0}, y expr=\thisrowno{1}] {dat_AT_TPM1f2x_bf_euc};
        \addplot[only marks, mark size = 0.5pt, color=blue] table[x expr=\thisrowno{0}, y expr=\thisrowno{1}] {dat_AT_TPM1f2x_adf_euc};
        \addplot[only marks, mark size = 0.5pt, color=orange] table[x expr=\thisrowno{0}, y expr=\thisrowno{1}] {dat_AT_TPM1f2x_haadf_euc};
        \addplot[only marks, mark size = 0.5pt, color=black] table[x expr=\thisrowno{0}, y expr=\thisrowno{1}] {dat_AT_TPM1f2x_3d_euc};
      \end{axis}
    \end{tikzpicture}
    \begin{tikzpicture}
      \begin{axis}[width = 0.49\textwidth,
                   xmin = 0, xmax = 1200,
                   ymin = 0, ymax = 0.5]
        \addplot[only marks, mark size = 0.5pt, color=green] table[x expr=\thisrowno{0}, y expr=\thisrowno{1}] {dat_AT_TPM1f2_bf_sup};
        \addlegendentry{BF}
        \addplot[only marks, mark size = 0.5pt, color=blue] table[x expr=\thisrowno{0}, y expr=\thisrowno{1}] {dat_AT_TPM1f2_adf_sup};
        \addlegendentry{ADF}
        \addplot[only marks, mark size = 0.5pt, color=orange] table[x expr=\thisrowno{0}, y expr=\thisrowno{1}] {dat_AT_TPM1f2_haadf_sup};
        \addlegendentry{HAADF}
        \addplot[only marks, mark size = 0.5pt, color=black] table[x expr=\thisrowno{0}, y expr=\thisrowno{1}] {dat_AT_TPM1f2_3d_sup};
        \addlegendentry{3D STEM}
        \addplot[only marks, mark size = 0.5pt, color=green] table[x expr=\thisrowno{0}, y expr=\thisrowno{1}] {dat_AT_TPM1f2x_bf_sup};
        \addplot[only marks, mark size = 0.5pt, color=blue] table[x expr=\thisrowno{0}, y expr=\thisrowno{1}] {dat_AT_TPM1f2x_adf_sup};
        \addplot[only marks, mark size = 0.5pt, color=orange] table[x expr=\thisrowno{0}, y expr=\thisrowno{1}] {dat_AT_TPM1f2x_haadf_sup};
        \addplot[only marks, mark size = 0.5pt, color=black] table[x expr=\thisrowno{0}, y expr=\thisrowno{1}] {dat_AT_TPM1f2x_3d_sup};
      \end{axis}
    \end{tikzpicture}\\
    \hspace*{-0.7em}
    \begin{tikzpicture}
      \begin{axis}[width = 0.49\textwidth,
                   legend pos = north west,
                   xmin = 0, xmax = 0.5,
                   ymin = 0, ymax = 0.5]
        \addplot[only marks, mark size = 0.5pt, color=green] table[x expr=\thisrowno{0}, y expr=\thisrowno{1}] {dat_AT_TPM1f2_io_bf_euc};
        \addlegendentry{BF}
        \addplot[only marks, mark size = 0.5pt, color=blue] table[x expr=\thisrowno{0}, y expr=\thisrowno{1}] {dat_AT_TPM1f2_io_adf_euc};
        \addlegendentry{ADF}
        \addplot[only marks, mark size = 0.5pt, color=orange] table[x expr=\thisrowno{0}, y expr=\thisrowno{1}] {dat_AT_TPM1f2_io_haadf_euc};
        \addlegendentry{HAADF}
        \addplot[only marks, mark size = 0.5pt, color=black] table[x expr=\thisrowno{0}, y expr=\thisrowno{1}] {dat_AT_TPM1f2_io_3d_euc};
        \addlegendentry{3D STEM}
        \addplot[only marks, mark size = 0.5pt, color=green] table[x expr=\thisrowno{0}, y expr=\thisrowno{1}] {dat_AT_TPM1f2x_io_bf_euc};
        \addplot[only marks, mark size = 0.5pt, color=blue] table[x expr=\thisrowno{0}, y expr=\thisrowno{1}] {dat_AT_TPM1f2x_io_adf_euc};
        \addplot[only marks, mark size = 0.5pt, color=orange] table[x expr=\thisrowno{0}, y expr=\thisrowno{1}] {dat_AT_TPM1f2x_io_haadf_euc};
        \addplot[only marks, mark size = 0.5pt, color=black] table[x expr=\thisrowno{0}, y expr=\thisrowno{1}] {dat_AT_TPM1f2x_io_3d_euc};
      \end{axis}
    \end{tikzpicture}\hspace*{0.7em}
    \begin{tikzpicture}
      \begin{axis}[width = 0.49\textwidth,
                   legend pos = north west,
                   xmin = 0, xmax = 0.5,
                   ymin = 0, ymax = 0.5]
        \addplot[only marks, mark size = 0.5pt, color=green] table[x expr=\thisrowno{0}, y expr=\thisrowno{1}] {dat_AT_TPM1f2_io_bf_sup};
        \addlegendentry{BF}
        \addplot[only marks, mark size = 0.5pt, color=blue] table[x expr=\thisrowno{0}, y expr=\thisrowno{1}] {dat_AT_TPM1f2_io_adf_sup};
        \addlegendentry{ADF}
        \addplot[only marks, mark size = 0.5pt, color=orange] table[x expr=\thisrowno{0}, y expr=\thisrowno{1}] {dat_AT_TPM1f2_io_haadf_sup};
        \addlegendentry{HAADF}
        \addplot[only marks, mark size = 0.5pt, color=black] table[x expr=\thisrowno{0}, y expr=\thisrowno{1}] {dat_AT_TPM1f2_io_3d_sup};
        \addlegendentry{3D STEM}
        \addplot[only marks, mark size = 0.5pt, color=green] table[x expr=\thisrowno{0}, y expr=\thisrowno{1}] {dat_AT_TPM1f2x_io_bf_sup};
        \addplot[only marks, mark size = 0.5pt, color=blue] table[x expr=\thisrowno{0}, y expr=\thisrowno{1}] {dat_AT_TPM1f2x_io_adf_sup};
        \addplot[only marks, mark size = 0.5pt, color=orange] table[x expr=\thisrowno{0}, y expr=\thisrowno{1}] {dat_AT_TPM1f2x_io_haadf_sup};
        \addplot[only marks, mark size = 0.5pt, color=black] table[x expr=\thisrowno{0}, y expr=\thisrowno{1}] {dat_AT_TPM1f2x_io_3d_sup};
      \end{axis}
    \end{tikzpicture}
  \end{center}
  \caption{Relative error of the final STEM image of the M1 specimen \cite{desanto04, blom18} calculated with LMA as compared to the result from the standard
           Multislice algorithm for STEM image simulation. Top row: relative error of the final STEM image for
           different values of $L$. Bottom row: relative error of the final STEM image in terms of the relative
           errors in the approximation of the STEM probes by linear combinations of
           the input waves. The relative errors are given in the euclidean norm
           (left column) and in the supremum norm (right column).
           The parameters are $u(x, y) = \varphi_n(x)\varphi_n(y)$ with $n=7$ and
           $(P_x, P_y) = (424, 400)$ with $X = Y = 1024$. The input wave lattice $\mathcal{I}_f$ is
           chosen as $\mathcal{I} = \mathcal{P} + \frac{1}{2}(p_x, p_y)$ with $f=2$.
           The simulation window size is given by $(l_x, l_y) = (84.54, 79.94)$ (in Angstrom) and the microscope settings are
           $\lambda = 0.0250793$ (corresponding to an accelerating
           voltage of $U = 200000$ Volt), $Z = 100$, $C_s = -2000$ and $\alpha_\text{max} = 0.026$.
           In order to reduce the computation time, only a $9\times 9$ pixel subsection of
           the final STEM image of size $P_x \times P_y$ was computed with a single frozen phonon iteration. This is
           because every data point in the above graphs corresponds to one full STEM image simulation. The integration domains used for the
           three different 2D STEM images are $[r_1, r_2] = [0, 15]$ (bright field, BF), $[r_1, r_2] = [16, 40]$ (annular dark field, ADF)
           and $[r_1, r_2] = [41, 200]$ (high angle annular dark field, HAADF) in mrad; cf. Step 3 in \cref{STEM_image_formation}.}\label{lma_inputoutputerror}
\end{figure}

In \cref{lma_inputoutputerror}, the error in the approximation of the STEM probe by a linear combination of trigonometric
polynomial input waves is compared to the error in the final output image for the M1
specimen \cite{desanto04, blom18} provided by Douglas Blom at University of South Carolina. Even for a relatively large relative error of -- for example -- 20\%
in the probe approximation, the relative euclidean error of the simulated STEM image is substantially
smaller than 20\% for the simulated BF, ADF and HAADF 2D STEM images as well as 3D STEM. Relative errors
measured in the supremum norm tend to be larger than relative errors in the euclidean norm, but are still
below 10\% for the BF and HAADF modes in this example.

\subsubsection{Partition of the STEM image}

Recall that each STEM probe position corresponds to a pixel in the final output image and that the
solution to the Schrödinger equation for a STEM probe is approximated by a linear combination of the
corresponding solutions for the input waves using the coefficients from \cref{lma_approximation_exit1} or
\cref{lma_approximation_exit2}. In this section, we assume that the probe and input wave lattices,
the simulation window pixel size, the approximation parameters $f$ and $L$ as well as the input wave
type have been chosen and that the approximation coefficients have been calculated.

One of the main differences of LMA and PRISM is that in LMA each STEM probe is approximated by a different
subset of all input waves, whereas in PRISM all input waves are used for the approximation of every
STEM probe. Since only a subset of $L$ input waves is needed at a time, it is possible to
calculate Multislice solutions for the input waves only when needed and to discard solutions that are not
needed anymore. The memory consumption and computation time can then be adjusted by choosing a suitable
order in which the output image pixels are computed.

The most straightforward approach minimizes the computation time: first, compute Multislice solutions for all input
wave positions $\mathcal{I}_f$; second, compute all linear combinations to calculate the output image pixels. This
requires storing $\abs{\mathcal{I}_f}$ images (one Multislice solution for each input wave) and $\abs{\mathcal{I}_f}$ executions of
the Multislice algorithm. A strategy that minimizes memory usage is to keep only those $L$ Multislice
solutions in memory that are required to form the linear combinations for the current pixel, discarding solutions
that are not needed anymore when continuing with the next pixel. This requires
storing only $L$ images and $k\in\N$ executions of the Multislice algorithm, where $\abs{\mathcal{I}_f}\le k \le \abs{\mathcal{I}_f}L$.
In general, much less than $\abs{\mathcal{I}_f}L$ executions will be necessary if the output image pixels are
computed in a sensible order. In LMA, the two sets of input wave positions used for the approximation
of neighboring probes generally have a large intersection due to the fact that the $L$ input waves
closest to the probe position are chosen for the approximation (cf. \cref{lrs_opt3} on \cpageref{lrs_opt3}).

Denote by $M\in\N$ the maximum number of Multislice solutions for the input waves that can
be kept in memory and assume that $M\ge L$. If $M$ is smaller than $\abs{\mathcal{I}_f}$, it is necessary
to determine an order in which the pixels of the output image are computed, where different orders may
require very different overall computation times. We again denote by $\mathcal{I}_{f,p,L}$ the $L$ input wave
positions in $\mathcal{I}_f$ closest to a probe position $p\in\mathcal{P}$, where $f\in\N$ is fixed, and extend
the notation to all subsets $S\subseteq\mathcal{P}$ by defining
\begin{equation*}
  \mathcal{I}_{f,S,L} := \bigcup_{p\in S} \mathcal{I}_{f,p,L}.
\end{equation*}

The idea is now to partition the set of the pixels $\mathcal{P}$ of the
STEM image into subsets $P_1,\ldots,P_l$
such that the number of Multislice solutions required for the computation of all pixels in a set $P_j$
satisfies $\abs{\mathcal{I}_{f,P_j,L}} \le M$
for all $j=1,\ldots,l$. If the pixels in $P_1$ are calculated first, then $P_2$, $P_3$ and so on, the total
number of executions of the Multislice algorithm is
\begin{equation*}
  T_{P_1,\ldots,P_l} := \sum_{j=1}^l \abs{\mathcal{I}_{f,P_j,L}} - \sum_{j=1}^{l-1} \abs{\mathcal{I}_{f,P_j,L} \cap \mathcal{I}_{f,P_{j+1},L}}.
\end{equation*}
The first sum gives the total number of Multislice solutions required by all of the sets $P_1,\ldots,P_l$ independently
of each other and the second sum corresponds to those Multislice solutions needed for the computation $P_j$ that can
be reused for $P_{j+1}$ and thus do not need to be recomputed.

Identifying a pixel $a(p_x\; 0) + b(0\;p_y)\in\mathcal{P}$ by its index $(a, b)$, we propose three
strategies to compute the Partition $\{P_1,\ldots,P_l\}$ of $\mathcal{P}$:
\begin{enumerate}
  \item \textbf{Row by row}: set $i=0$, $(a_0, b_0) := (0, 0)$ and $P_1 := {(a_0, b_0)}$. While $\abs{\mathcal{I}_{f,P_1,L}} \le M$,
        set
        \begin{equation*}
          (a_{i+1}, b_{i+1}) := \begin{cases} (a_i+1, b_i), & \text{if } a_i < P_x-1, \\
                                              (0, b_i+1), & \text{otherwise} \end{cases}
        \end{equation*}
        and $P_1 := P_1 \cup (a_{i+1}, b_{i+1})$. Continue with $P_2, \ldots, P_l$ in the same fashion.
  \item \textbf{Rectangles}: partition the set $\mathcal{P}$ into rectangular sets $P_1,\ldots,P_l$
        of the same size that are as large as possible while satisfying the condition
        $\abs{\mathcal{I}_{f,P_j,L}} \le M$ at the same time; fill any remaining space with smaller
        rectangles. The order of the rectangles in the partition can be row-major or column-major.
  \item \textbf{Greedy algorithm}: choose a random starting point $(a_0, b_0)$ and set $P_0 = {(a_0, b_0)}$.
        New pixels are added to $P_0$ by choosing the neighboring pixel $(a_{i+1}, b_{i+1})$ to the
        previously added pixel $(a_i, b_i)$ with the lowest cost, i.e. such that
        \begin{equation*}
          \abs{\mathcal{I}_{f, P_0, L} \cup \mathcal{I}_{f, (a_{i+1}, b_{i+1}), L}}
        \end{equation*}
        is minimal, where $a_i-1\le a_{i+1} \le a_i+1$ and $b_i-1\le b_{i+1}\le b_i+1$.
        Continue for as long as possible while $\abs{\mathcal{I}_{f, P_0, L}} \le M$ and then
        continue with $P_1,P_2,\ldots,P_l$ in the same way, starting with the last pixel that
        could not be added anymore to the previous set.
\end{enumerate}

\begin{figure}
  \begin{minipage}{0.32\textwidth}
    \includegraphics[width=\textwidth, keepaspectratio]{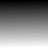}\\
    \vspace*{-2em}
    \begin{flushright} $T_{P_1,\ldots,P_l} = 1642$ ($l=60$)\end{flushright}
  \end{minipage}\hfill
  \begin{minipage}{0.32\textwidth}
    \includegraphics[width=\textwidth, keepaspectratio]{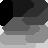}\\
    \vspace*{-2em}
    \begin{flushright} $T_{P_1,\ldots,P_l} = 1856$ ($l=12$)\end{flushright}
  \end{minipage}\hfill
  \begin{minipage}{0.32\textwidth}
    \includegraphics[width=\textwidth, keepaspectratio]{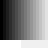}\\
    \vspace*{-2em}
    \begin{flushright} $T_{P_1,\ldots,P_l} = 1276$ ($l=18$)\end{flushright}
  \end{minipage}\\[1em]
  
  \begin{minipage}{0.32\textwidth}
    \includegraphics[width=\textwidth, keepaspectratio]{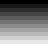}\\
    \vspace*{-2em}
    \begin{flushright} $T_{P_1,\ldots,P_l} = 953$ ($l=12$)\end{flushright}
  \end{minipage}\hfill
  \begin{minipage}{0.32\textwidth}
    \includegraphics[width=\textwidth, keepaspectratio]{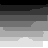}\\
    \vspace*{-2em}
    \begin{flushright} $T_{P_1,\ldots,P_l} = 995$ ($l=10$)\end{flushright}
  \end{minipage}\hfill
  \begin{minipage}{0.32\textwidth}
    \includegraphics[width=\textwidth, keepaspectratio]{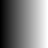}\\
    \vspace*{-2em}
    \begin{flushright} $T_{P_1,\ldots,P_l} = 953$ ($l=48$)\end{flushright}
  \end{minipage}
  \caption{Two examples for partitions computed with the strategies row by row (left column),
           greedy (center column) and rectangle (right column). Different gray values correspond to
           different sets of the partition, where $P_1$ consists of all black pixels and $P_l$ of
           all white pixels. The minimum number of
           Multislice calls that can be achieved, $T_\mathcal{P} = \abs{\mathcal{I}_{f,\mathcal{P},L}}$, is 1077 (top row)
           respectively 953 (bottom row). Both experiments only differ in the target approximation
           error of the STEM probe, which is $0.01$ for the experiment shown in the top row and
           $0.03$ for the experiment in the bottom row, resulting in $L=75$ (approximation error $0.01$)
           and $L=51$ (approximation error $0.03$. The other
           parameters are $u=\hat\psi_0^\text{init}$ and $\mathcal{I} = \mathcal{P} + \frac{1}{2}(p_x, p_y)$
           with $f=2$ and $X = Y = 512$ as well as $P_x = P_y = 158$, where only a subsection of the STEM image of
           size $48 \times 48$ was to be computed. The simulation window size is $l_x = l_y = 31.208$ (in Angstrom)
           and the microscope settings are $\lambda = 0.0250793$ (corresponding to an accelerating
           voltage of $U = 200000$ Volt), $Z = 100$, $C_s = -2000$ and $\alpha_\text{max} = 0.026$.}\label{lmastrategy_examples}
\end{figure}

\cref{lmastrategy_examples} shows two examples for partitions computed with each of these
three strategies. It is interesting to note that this approach doesn't necessarily increase
the computation time: the partitions computed with the row by row and rectangle
strategies in the second row of \cref{lmastrategy_examples} result in exactly the same
number of Multislice calls as the trivial strategy $P_1 = \mathcal{P}$, but require significantly
less computer memory. The partition computed with the rectangle strategy only requires storing
about 1/4th of the Multislice results at any time as compared to the trivial strategy,
reducing the required computer memory by 75\% without increasing the computation time.

\subsubsection{Comparison of different input wave types}

\begin{figure}
  \begin{center}
  \includegraphics[width=0.29\textwidth, keepaspectratio]{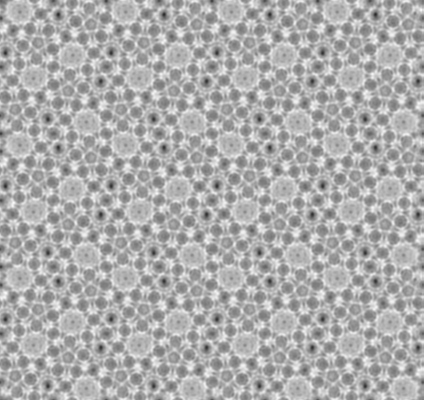}
  \includegraphics[width=0.29\textwidth, keepaspectratio]{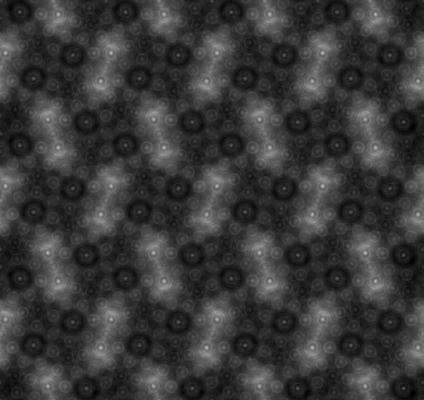}
  \includegraphics[trim={6cm 6.345cm 6cm 5cm},clip,width=0.29\textwidth, keepaspectratio]{dat_LMAE_M12_HAADF}\\[0.2em]
  
  \includegraphics[width=0.29\textwidth, keepaspectratio]{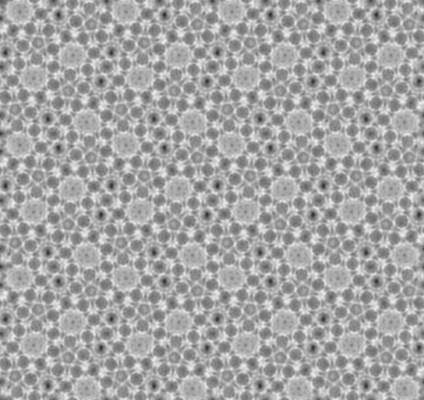}
  \includegraphics[width=0.29\textwidth, keepaspectratio]{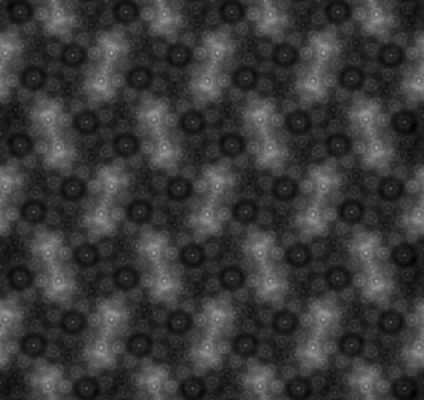}
  \includegraphics[trim={6cm 6.32cm 6cm 5cm},clip,width=0.29\textwidth, keepaspectratio]{dat_LMAE_M13_HAADF}\\[0.2em]
  
  \includegraphics[width=0.29\textwidth, keepaspectratio]{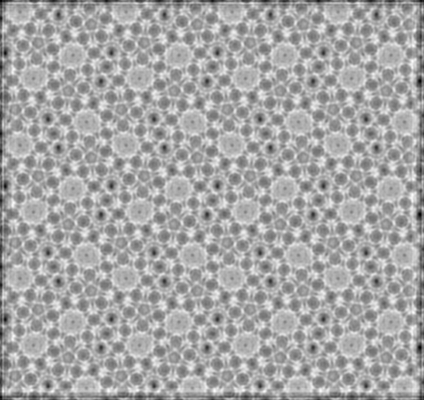}
  \includegraphics[width=0.29\textwidth, keepaspectratio]{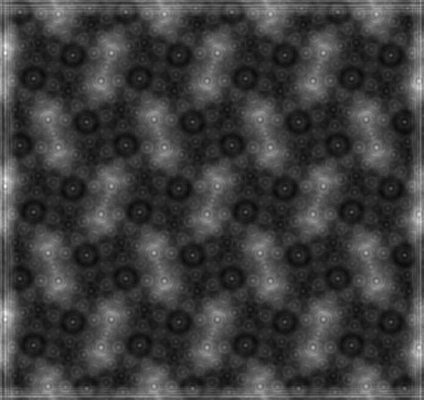}
  \includegraphics[trim={6cm 6.32cm 6cm 5cm},clip,width=0.29\textwidth, keepaspectratio]{dat_LMAE_M14_HAADF}\\[0.2em]
  
  \includegraphics[width=0.29\textwidth, keepaspectratio]{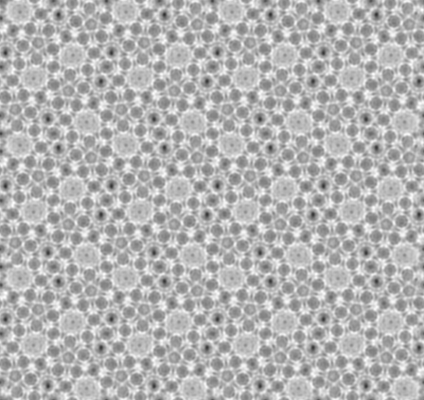}
  \includegraphics[width=0.29\textwidth, keepaspectratio]{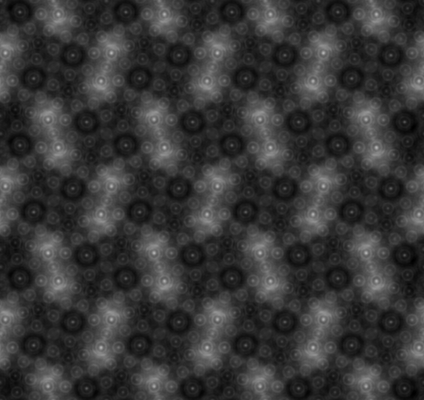}
  \includegraphics[trim={6cm 6.32cm 6cm 5cm},clip,width=0.29\textwidth, keepaspectratio]{dat_LMAE_M15_HAADF}
  \end{center}
  \caption{Simulated images of the M1 specimen with the LMA algorithm and four different configurations.
           Top row: $u=\hat\psi_0^\text{init}$ with $\mathcal{I} = \mathcal{P}$, $f=1$ and $L=1$.
           Second row: $u=\hat\psi_0^\text{init}$ with $\mathcal{I} = \mathcal{P} + \frac{1}{2}(p_x, p_y)$, $f=2$ and $L=27$.
           Third row: $u(x,y)=\varphi_n(x)\varphi_n(y)$ with $\mathcal{I} = \mathcal{P} + \frac{1}{2}(p_x, p_y)$, $f=2$ and $L=283$.
           Bottom row: $u=g_\sigma$ with $\sigma=0.482$ Angstrom, $\mathcal{I} = \mathcal{P} + \frac{1}{2}(p_x, p_y)$, $f=2$ and $L=379$.
           Left column: BF images ($[r_1, r_2] = [0, 15]$). Middle column: HAADF images ($[r_1, r_2] = [41, 200]$).
           Right column: subsection of the HAADF images. The relative probe approximation errors in the
           supremum norm are 0\% for the top row and about 10\% for the other rows.}\label{lmaexample_figure}
\end{figure}

\cref{lmaexample_figure} shows the simulation results for the M1 specimen \cite{desanto04, blom18} with
four different configurations of LMA. As expected from the results in \cref{lma_inputoutputerror}, the
results look almost identical to the reference images in the top row, even though an relative approximation
error of about 10\% in the supremum norm was used in the simulation to keep $L$ small and the simulation
fast. The simulation of the images in the top row took about 24 hours, whereas the second row took 8 hours,
the third row 12 hours and the fourth row 13 hours.

\begin{figure}
  \includegraphics[width=0.32\textwidth, keepaspectratio]{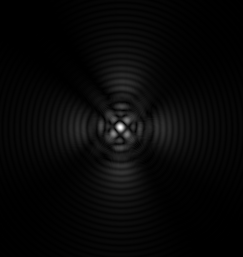}
  \includegraphics[width=0.32\textwidth, keepaspectratio]{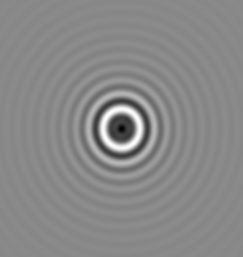}
  \includegraphics[width=0.32\textwidth, keepaspectratio]{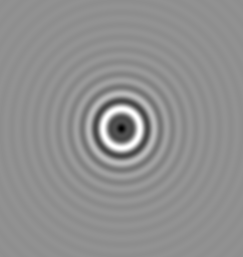}\\[-0.75em]
  
  \includegraphics[width=0.32\textwidth, keepaspectratio]{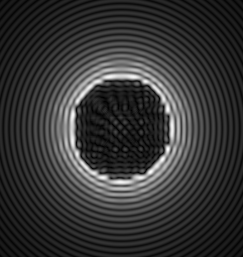}
  \includegraphics[width=0.32\textwidth, keepaspectratio]{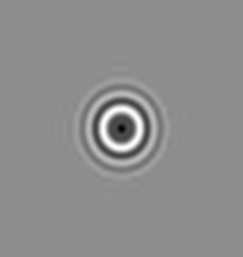}
  \includegraphics[width=0.32\textwidth, keepaspectratio]{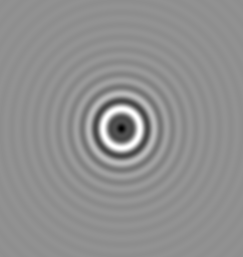}\\[-0.75em]
  
  \includegraphics[width=0.32\textwidth, keepaspectratio]{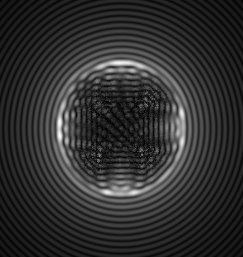}
  \includegraphics[width=0.32\textwidth, keepaspectratio]{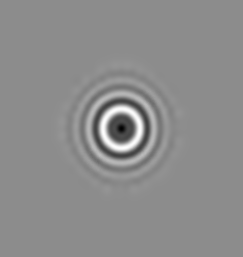}
  \includegraphics[width=0.32\textwidth, keepaspectratio]{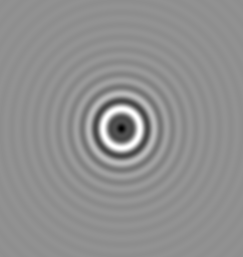}
  \caption{Approximations of the STEM probes corresponding to the second, third and fourth row of
           \cref{lmaexample_figure}. Top row: $u=\hat\psi_0^\text{init}$. Middle row: $u(x,y)=\varphi_n(x)\varphi_n(y)$.
           Bottom row: $u=g_\sigma$ with $\sigma=0.482$ Angstrom. Left column: relative errors in the supremum norm, where
           white equals 10\% and black equals 0\%. Middle column: Imaginary part of the resulting approximation of the
           STEM probe. Right column: imaginary part of the STEM probe.}\label{lmaexample_approx}
\end{figure}

It is interesting to see the very different probe approximations that were used in these simulations
as shown in \cref{lmaexample_approx}. In particular, well localized functions such as Gaussians or
trigonometric polynomials only approximate the probes near their center. This will not yield a significant
error in the output image as was already established by PRISM, where the results from the Multislice computations
are cropped by a factor of $f^2$ around the probe, which has a similar effect.

\subsubsection{Recomputing local changes}

One of the strong suits of the approach presented here is that the small support of the input
waves in real space allows for a fast recomputation if local changes are made to the specimen such as
changing a single atomic column. The reason for this is that we only need to redo those Multislice
computations that are needed for the approximation of the probe positions close to the location of
the local changes. If $\mathcal{P}'\subseteq\mathcal{P}$ is the set of pixels of the STEM image that
need to be recomputed, then the Multislice algorithm needs to be performed again 
only for those input wave positions in $\mathcal{I}_{f,\mathcal{P}',L}$ that are close enough to the
changes in the specimen to have the result of the Multislice computation changed.

In practice, this reduces the computation time twice: first, because we only need to recompute pixels
within a small area around the changes; second, because we similarly only need to redo Multislice
computations for input wave positions within a small area around the changes. While it is possible
to compute individual pixels of the output image with PRISM, this always requires performing the Multislice
algorithm for all input waves on the entire simulation window. MULTEM on the other hand requires performing
the Multislice algorithm once for every probe position, providing no additional speedup.

\begin{figure}
  \begin{tikzpicture}
    \node at (0, 0) {\includegraphics[trim={0 6cm 8cm 1cm},clip,width=0.31\textwidth, keepaspectratio]{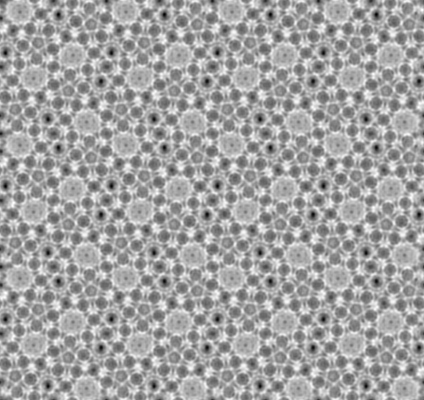}};
    \draw[->,color = red, line width = 1pt] (0.3, 0.5) -> (-0.2, 0.2);
  \end{tikzpicture}
  \begin{tikzpicture}
    \node at (0, 0) {\includegraphics[trim={0 6cm 8cm 1cm},clip,width=0.31\textwidth, keepaspectratio]{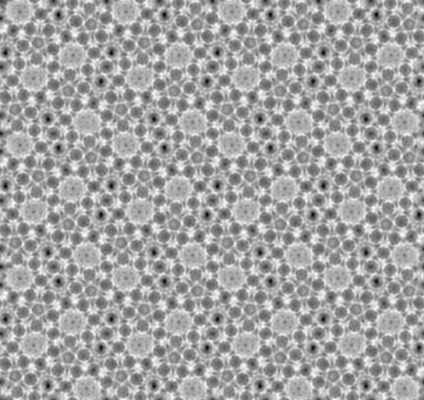}};
    \draw[->,color = red, line width = 1pt] (0.3, 0.5) -> (-0.2, 0.2);
  \end{tikzpicture}
  \begin{tikzpicture}[blend group=multiply]
    \node at (0, 0) {\includegraphics[trim={0 6cm 8cm 1cm},clip,width=0.31\textwidth, keepaspectratio]{dat_Recomp_BF4}};
    \node at (0, 0) {\includegraphics[trim={0 6cm 8cm 1cm},clip,width=0.31\textwidth, keepaspectratio]{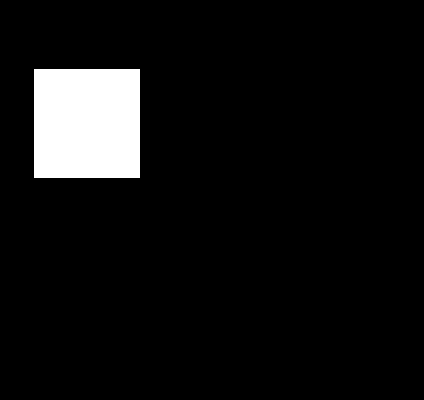}};
  \end{tikzpicture}
  
  \begin{tikzpicture}
    \node at (0, 0) {\includegraphics[trim={0 6cm 8cm 1cm},clip,width=0.31\textwidth, keepaspectratio]{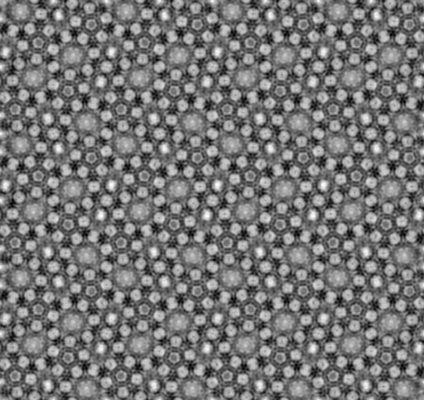}};
    \draw[->,color = red, line width = 1pt] (0.3, 0.5) -> (-0.2, 0.2);
  \end{tikzpicture}
  \begin{tikzpicture}
    \node at (0, 0) {\includegraphics[trim={0 6cm 8cm 1cm},clip,width=0.31\textwidth, keepaspectratio]{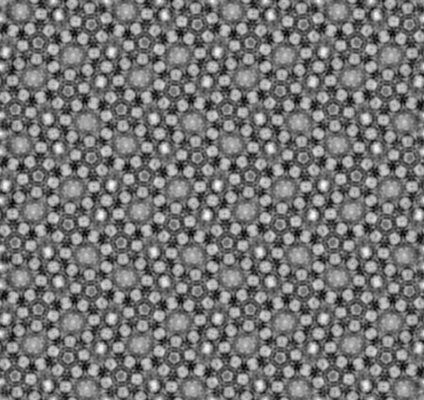}};
    \draw[->,color = red, line width = 1pt] (0.3, 0.5) -> (-0.2, 0.2);
  \end{tikzpicture}
  \begin{tikzpicture}[blend group=multiply]
    \node at (0, 0) {\includegraphics[trim={0 6cm 8cm 1cm},clip,width=0.31\textwidth, keepaspectratio]{dat_Recomp_ADF4}};
    \node at (0, 0) {\includegraphics[trim={0 6cm 8cm 1cm},clip,width=0.31\textwidth, keepaspectratio]{dat_Recomp_probe_update_positions4}};
  \end{tikzpicture}
  
  \begin{tikzpicture}
    \node at (0, 0) {\includegraphics[trim={0 6cm 8cm 1cm},clip,width=0.31\textwidth, keepaspectratio]{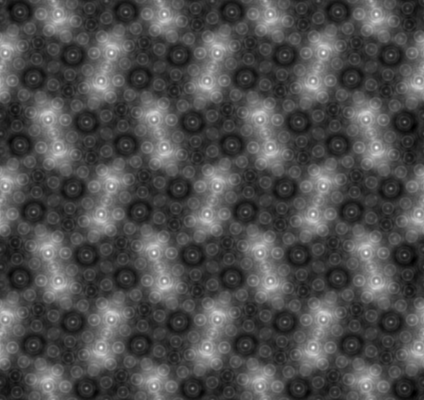}};
    \draw[->,color = red, line width = 1pt] (0.3, 0.5) -> (-0.2, 0.2);
  \end{tikzpicture}
  \begin{tikzpicture}
    \node at (0, 0) {\includegraphics[trim={0 6cm 8cm 1cm},clip,width=0.31\textwidth, keepaspectratio]{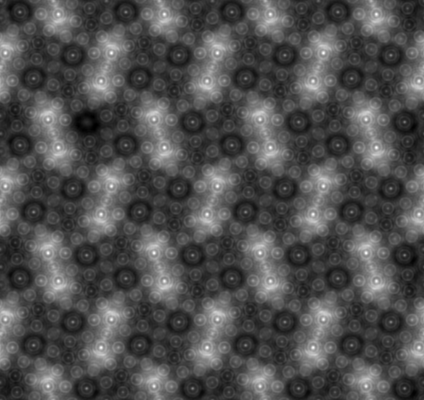}};
    \draw[->,color = red, line width = 1pt] (0.3, 0.5) -> (-0.2, 0.2);
  \end{tikzpicture}
  \begin{tikzpicture}[blend group=multiply]
    \node at (0, 0) {\includegraphics[trim={0 6cm 8cm 1cm},clip,width=0.31\textwidth, keepaspectratio]{dat_Recomp_HAADF4}};
    \node at (0, 0) {\includegraphics[trim={0 6cm 8cm 1cm},clip,width=0.31\textwidth, keepaspectratio]{dat_Recomp_probe_update_positions4}};
  \end{tikzpicture}
  \caption{STEM images of the M1 specimen simulated with LMA. Top row: bright field ($[r_1, r_2] = [0, 15]$);
           middle row: annular dark field ($[r_1, r_2] = [16, 40]$); bottom row: high angle annular dark field
           ($[r_1, r_2] = [41, 200]$). The images in the left and middle column only differ by changes
           to the atomic column indicated by the red arrow. For the images in the left column, the atomic column
           consists of alternating Te and O atoms, whereas all Te atoms have been replaced by O atoms for the
           image in the middle column. The images in the right column show those pixels that had to be
           recomputed, i.e. the images in the middle column are a combination of the images in the left
           and right column. Shown here are cropped subsections of the specimen for a better visibility
           of the differences between the images.}\label{recompute_M1}
\end{figure}

\cref{recompute_M1} shows an example where two configurations of the M1 specimen \cite{desanto04, blom18} are simulated that differ only
in one atomic column. In this example, the simulation window size is $1024 \times 1024$ pixels with a
physical size of $84.536 \times 79.941$ Angstrom. The size of the probe lattice is $(P_x, P_y) = (424, 400)$ and
$\mathcal{I} = \mathcal{P} + \frac{1}{2}(p_x, p_y)$ with $f=2$. The maximum relative approximation error in
the supremum norm in this example is set to $2\%$, which results in $L=59$ when computing the
approximation coefficients with linear least squares.

For the images shown in the left column of \cref{recompute_M1}, this results in
$P_xP_y = 169600$ probe positions and an overall number of $\abs{\mathcal{I}_f} = 212\cdot 200 = 42400$
times that the Multislice algorithm needs to be perfomed. If one was then interested in changes in the
intensity if only the atomic column depicted by the red arrow in \cref{recompute_M1} was changed,
it would be a waste to repeat the entire computation.
Instead, using the approach outlined above, only 11554 probe positions need to be recomputed,
which requires the recomputation of 3825 Multislice solutions. In our experiment, this reduced
the computation time from 150 minutes to 14 minutes for every image after the initial one.

\paragraph{Source code:} the source code that implements LMA and has been used for all experiments
in this paper is made freely available at \url{https://github.com/CDoberstein/LMA}.

 \vskip .1in
 \noindent
\textbf{Acknowledgments:} The authors thank Thomas Vogt, Douglas Blom and Michael Matthews for helpful discussions 
on the research in this paper.

\bibliographystyle{plain}
\bibliography{main}

\end{document}